\documentclass{USG_nlaa}
\usepackage{graphicx}
\usepackage{amsmath}
\usepackage{mathtools}
\usepackage{algorithm}
\numberwithin{algorithm}{section}
\usetikzlibrary{decorations.pathreplacing}

\makeatletter
\algrenewcommand\ALG@beginalgorithmic{\normalfont\fontsize{8.5}{10.5}\selectfont}
\def\oddhead@titlepage@info{}
\def\evenhead@titlepage@info{}
\makeatother

\articletype{Research Article}

\historydates{}
\journal{Numerical Linear Algebra with Applications}

\hypersetup{pdftitle={Conjugation-preserving polynomial and rational interpolation for limited-memory matrix-function actions},
  pdfauthor={Tianqi Zhang and Fei Xue},
  pdfsubject={Limited-memory approximation of matrix-function actions},
  colorlinks=true,allcolors=black}

\begin{document}

\title{Conjugation-preserving polynomial and rational interpolation for limited-memory matrix-function actions}

\author[1]{Tianqi Zhang}

\author[2]{Fei Xue}

\authormark{T. Zhang and F. Xue}

\address[1]{\orgdiv{School of Data Sciences}, \orgname{Zhejiang University of Finance and Economics}, \orgaddress{\city{Hangzhou}, \state{Zhejiang}, \country{China}}}

\address[2]{\orgdiv{School of Mathematical and Statistical Sciences}, \orgname{Clemson University}, \orgaddress{\city{Clemson},\state{SC}, \country{USA}}}

\corres{Fei Xue, School of Mathematical and Statistical Sciences, Clemson University, O-110 Martin Hall, Box 340975, Clemson, SC 29634.\ \email{fxue@clemson.edu}}

\fundingInfo{U.S. National Science Foundation under grant DMS-2111496.}

\abstract[Summary]{Approximating $f(A)b$ for a large nonsymmetric matrix can require substantial storage if a long Krylov basis is retained. We develop polynomial and rational interpolation schemes that use a fixed number of large recurrence vectors and preserve conjugation symmetry. Starting from a prescribed scalar approximation region, the schemes organize real nodes and complex-conjugate pairs into Newton-like bases with real interpolation coefficients. The polynomial and single-real-pole recurrences use real arithmetic in the matrix-action stage; the multi-pole construction also exploits conjugate pole pairs. A conjugation-preserving AAA procedure supplies interpolation data for the latter. We establish the interpolation and basis properties and compare the resulting realizations with Leja, Taylor, and Krylov methods on matrices from differential equations, model reduction, and networks. The experiments demonstrate arithmetic savings relative to complex-arithmetic counterparts and identify the influence of shifted solves and orthogonalization requirements. The methods are intended for memory-constrained computations in which useful approximation-region information is available.}

\keywords{matrix-function actions; polynomial interpolation; rational interpolation; short recurrences; real arithmetic; AAA algorithm}

\jnlcitation{\cname{
\author{T. Zhang}, and
\author{F. Xue}} (\cyear{2026}),
\ctitle{Conjugation-preserving polynomial and rational interpolation for limited-memory matrix-function actions}, \cjournal{Numer. Linear Algebra Appl.}, \cvol{}.}

\maketitle

\section{Introduction}\label{sec:introduction}
Computing the action of a matrix function on a vector, $f(A)b$, is a basic task in scientific computing; see \cite{Gue13b,Higham08book} and the references therein. We consider a large, sparse, generally nonsymmetric matrix $A\in\mathbb R^{n\times n}$ and a nonzero vector $b\in\mathbb R^n$. Our aim is to approximate this action without retaining a long basis of large vectors, while taking advantage of the real input data even when complex interpolation points and poles are useful.

The methods start from a compact scalar approximation region $D\subset\mathbb C$ that is symmetric with respect to the real axis. In the analytic setting, $f$ is holomorphic on a neighborhood of $D$ and satisfies $f(\overline z)=\overline{f(z)}$. We construct polynomial or rational interpolants on $D$ and evaluate their matrix actions by short recurrences. Here, \emph{limited memory} refers to a fixed number of length-$n$ recurrence vectors, independent of the approximation degree. Scalar interpolation data and the storage required by shifted linear solvers are additional costs.

A standard alternative constructs an orthonormal basis $V_m$ for a polynomial or rational Krylov subspace,
\[
 \mathcal K_m(A,b)=\operatorname{span}\{b,Ab,\ldots,A^{m-1}b\},\qquad
 \mathcal Q_m(A,b)=\prod_{\ell=1}^{m-1}(A-\xi_\ell I)^{-1}\mathcal K_m(A,b),
\]
and forms $f_m=V_mf(H_m)e_1\|b\|_2$, where $H_m=V_m^TAV_m$ for a real orthonormal basis. These methods are well established, but full orthogonalization costs $\mathcal O(nm^2)$ operations and retaining the basis costs $\mathcal O(nm)$ storage. Restarting reduces the storage requirement \cite{AfanasjewEtal08,EiermannErnst06,EiermannEtal11,FrommerEtal14a,FrommerEtal14b}. Randomized sketching can reduce orthogonalization costs \cite{GuttelSimunec24,CortinovisEtal24,GuttelSchweitzer23,BurkeEtal24}, while two-pass realizations regenerate the basis instead of storing it \cite{GuttelSchweitzer23,PalittaSchweitzerSimoncini25}. The balance between arithmetic, storage, and communication is particularly important in large-scale computations \cite{DemmelEtal12}.

Interpolation provides a complementary route. Existing limited-memory methods include polynomial and rational Leja schemes \cite{CaliariEtal04,CaliariEtal16,RostamiXue18}, truncated Taylor expansions \cite{AlMohy2011,HighamKandolf17,AlMohy18}, rational Lanczos \cite{PalittaEtal22}, and the poles-and-interpolation-nodes algorithm (PAIN) \cite[section~5.4.2]{GuttelThesis}. PAIN is the closest predecessor to our multi-pole construction: both evaluate the same final rational interpolant when given the same admissible nodes and poles. The contribution here is a systematic \emph{conjugation-preserving realization} of polynomial, single-pole, and multi-pole interpolation. We group nonreal data into conjugate pairs, construct Newton-like bases with real coefficients, and derive short recurrences for the matrix action. A symmetry-preserving form of AAA supplies the multi-pole interpolation data. The resulting matrix operations are real for the polynomial and single-real-pole schemes, and exploit conjugate structure for the multi-pole scheme. This use of conjugation follows the same broad principle as real-arithmetic treatments of conjugate eigenpairs \cite[section~7.5.4]{golub2013matrix}\cite{parlett1987complex}.

The region $D$ is an input, not an output of the matrix recurrence. This is useful when approximation-region information is available from the application or can be reused for several right-hand sides. In return for this information, the methods avoid basis orthogonalization and do not require an orthogonalization-window parameter. They are not intended to replace adaptive Krylov methods when no useful region is available. Section~\ref{sec:experiment} compares these different requirements explicitly, including the matrix-dependent windows observed in a two-pass truncated-Arnoldi implementation.

The link between scalar approximation and the matrix action is given by the Crouzeix--Palencia theorem \cite{CrouzeixPalencia17}. Define
\[
 W(A)=\{x^*Ax:x\in\mathbb C^n,\ \|x\|_2=1\}.
\]
If $f$ is holomorphic on a neighborhood of $W(A)$ and $q_m$ is a polynomial, or a rational function with no poles on that neighborhood, then
\begin{equation}\label{crouzeixbound}
 \|f(A)-q_m(A)\|_2\leq(1+\sqrt2)\sup_{z\in W(A)}|f(z)-q_m(z)|.
\end{equation}
Consequently, with $\|h\|_D=\sup_{z\in D}|h(z)|$,
\[
 \|f(A)b-q_m(A)b\|_2\leq(1+\sqrt2)\|f-q_m\|_D\,\|b\|_2
 \quad\text{if }W(A)\subset D.
\]
This bound explains the value of a suitable enclosing region. A region inferred from a few eigenvalues need not enclose $W(A)$, and a small residual on a finite sample set is not a bound for the continuous norm. Our experiments therefore report observed errors alongside the prescribed regions, stopping tests, and action-stage costs. The singular graph-Laplacian examples use a separate interpretation of $-z\log z$ at zero, described in section~\ref{subsec:AlogAb}.

Section~\ref{sec:review} reviews related limited-memory methods. Section~\ref{sec:asymfac} explains the role of conjugation symmetry. Sections~\ref{sec:poly}--\ref{sec:lmrksm} develop the polynomial, single-pole, and multi-pole recurrences and establish their interpolation properties. Section~\ref{sec:experiment} presents the numerical comparisons and discusses their implications, and section~\ref{sec:conclusion} concludes.

\section{Related limited-memory methods}\label{sec:review} We review several limited-memory polynomial and rational approximations to $f(A)b$. The review clarifies their different memory, arithmetic, and parameter-selection requirements. 

\subsection{Real Leja polynomial interpolation}\label{subsec:realleja} Let $K\subset\mathbb C$ be a regular compact continuum of positive capacity whose complement in the Riemann sphere is simply connected (for example, a filled ellipse), $\mathbb D\subset\mathbb C$ be the closed unit disk centered at the origin, and ${\overline{\mathbb{C}}} = \mathbb{C} \cup \{\infty\}$. Let $\phi$ be the Riemann mapping that maps ${\overline{\mathbb{C}}} \setminus K$ conformally onto ${\overline{\mathbb{C}}} \setminus\mathbb D$, such that $\phi(\infty) = \infty$ and $\lim_{z\rightarrow \infty}\phi(z)/z = \frac{1}{\gamma}$, where $\gamma = \mathrm{cap}(K)$; { see \cite[Theorem 3.14]{markushevich1967} and \cite[sections~1.1 \& 2.2]{Gaier1987}}. For any $r > \gamma$, define $\Gamma_r = \{z: |\phi(z)| = \frac{r}{\gamma}\}$ and let $K_r$ be the bounded domain with boundary $\Gamma_r$. We set $K_\gamma=K$, and $K_{r_1} \subset K_{r_2}$ if $r_1 \leq r_2$. Let $R>\gamma$ be the supremum of the parameters $r$ for which $f$ is holomorphic in the interior of $K_r$, allowing $R=\infty$ for an entire function. Let $z_0\in K$ be fixed. The Leja points {\cite{Leja1957,Reichel1990}} are defined recursively as
\begin{eqnarray}\label{eq:lejapts}
    z_{m} = {\mathrm{argmax}_{z\in K}} \textstyle \prod_{j=0}^{m-1}|z-z_j|.\quad(m=1,2,\ldots)
\end{eqnarray}
The sequence of polynomials $\{p_m\}$ interpolating $f(z)$ at the Leja points $z_0,\ldots,z_m \in K = K_\gamma$ converges at the maximal geometric rate on $K$, namely $\limsup_{m\rightarrow\infty}\|f-p_m\|_K^{1/m}=\gamma/R${\cite{Reichel1990,Walsh1935}}. Let $K=K_\gamma$ be the ellipse centered at $d \in \mathbb{R}$, with semiaxes $a$ and $b$, and foci $d\pm c$, such that the capacity of $K$ is $\gamma = \frac{a+b}{2}$. For the horizontal-major-axis case, extend the notation $K_r$ to the confocal family down to its degenerate member $K_{c/2}=[d-c,d+c]$. It is shown that a sequence of polynomials $\{p_m\}$ converging maximally to $f(z)$ on $[d-c,d+c]$ also converges to $f(z)$ on $K_r$ for some $r \geq \frac{c}{2}$ \cite[section 3, Theorems 22 \& 23]{NovatiThesis}, and the asymptotic convergence factor is $\frac{r}{R}$ depending on $K_r$.

The real Leja polynomial interpolation \cite{CaliariEtal04} needs a rectangle encompassing the spectrum $\sigma(A)$. It constructs the ellipse with the smallest capacity that circumscribes this rectangle. The real Leja points $\{z_j\}_{j=0}^m = \left\{d+\frac{c}{2} \xi_j\right\}_{j=0}^m$ are defined on the focal interval $K_{c/2}=[d-c,d+c]$, and $\{\xi_j\}$ are precomputed Leja points on $[-2,2]$. The resulting Newton interpolation can be evaluated by short recurrences. {We refer to \cite{FischerReichel1989,Reichel1990} for the importance of ordering of interpolation points for the conditioning of Newton interpolation, and \cite[section~5.3]{Higham2002} for a general error analysis of divided differences.} 
If the focal interval is parallel to the imaginary axis, let $z_0=d \in \mathbb{R}$ be the center of the ellipse, then use \eqref{eq:lejapts} to obtain $z_m$, and define $z_{m+1} = \overline{z}_m$ to obtain conjugation-invariant Leja points. With a careful construction, one can avoid complex arithmetic to construct Newton interpolations of even degrees. \smallskip

\paragraph{Function-specific ingredients.}
The real Leja implementation used here targets the exponential, although polynomial interpolation is not restricted to that function. Two ingredients are useful in this setting: accurate exponential divided differences \cite{McCurdy1984}, and substepping through
\[
 \exp(hA)b=\bigl[\exp((h/s)A)\bigr]^sb.
\]
The integer $s$ reduces the difficulty of each action. Other functions require their own divided-difference evaluation and, where available, an analogous scaling identity. If $p_m$ approximates $e^z$, the corresponding action approximation is $p_m(hA)b$.

\subsection{Truncated Taylor series}
Taylor expansions provide limited-memory action methods for exponential, trigonometric, and hyperbolic functions. For the exponential, the basic polynomial is combined with scaling and repeated application \cite{AlMohy2011}. 

The action is approximated as
\[
 \exp(A)b=\bigl(\exp(A/s)\bigr)^s b
 \approx p_m(A/s)^s b,\qquad p_m(z)=\sum_{j=0}^m z^j/j!.
\]
The degree $m$ and positive integer $s$ are chosen together using norm-based backward-error bounds, including estimates involving powers of $A$. The resulting action is evaluated by matrix-vector products and repeated application to a vector, without forming the full matrix exponential. This distinction matters for nonnormal matrices: the spectral radius alone does not control the truncation or backward error.

\subsection{Single-pole rational Chebyshev and Leja methods} Approximations to $f(z) = e^z$ on $(-\infty,0]$ can be achieved by rational functions with a {repeated pole \cite[Chapter 25]{Trefethen19book}}. The single-pole {rational Chebyshev method \cite{Guttel.RatChebExp.2010}} determines a scalar $\alpha>0$ and defines the new variable $w = \frac{\alpha + z}{\alpha -z} \in (-1,1]$, i.e., $z = \alpha \frac{w-1}{w+1}$, so that $e^z = \mathrm{exp}\big(\alpha \frac{w-1}{w+1}\big):= g(w)$. {For a smaller tolerance, a larger value $\alpha \in [0.57,27.29]$ and a higher degree $m$ (up to $36$) are adopted from a table of prescribed data, so that $g(w) \approx \sum_{j=0}^m c_j T_j(w)$ within the desired tolerance. The coefficients $c_j$ are computed by fast Fourier transforms \cite[Chapter 3]{Trefethen19book}. Replacing $z$ and $w$ with $A$ and $M=(\alpha I - A)^{-1}(\alpha I + A)$ respectively, we have $\exp(A)b = g(M)b \approx \sum_{j=0}^m c_j T_j(M)b$. The three-term recurrence of Chebyshev polynomials is used to evaluate $T_j(M)b$ of an increasing degree $j$. The computation of $\sum_{j=0}^m c_j T_j(M)b$ needs $m$ evaluations} of a vector of the form $Mu = (\alpha I - A)^{-1}(\alpha I + A)u$, achieved by solving linear systems with the fixed coefficient matrix $\alpha I - A$. For a real symmetric negative-semidefinite matrix, the scalar approximation error directly bounds the matrix error.

Similarly, the single-pole rational Leja method \cite{RostamiXue18} constructs polynomial interpolation $p_m(w)$ in Newton basis at the Leja points on $[-1,1]$ to approximate $g(w)=e^z$. {It adopts a fixed $\alpha = 50$ and the highest degree $m=45$, determined experimentally by the performance on several nonsymmetric matrices \cite[section 3.3.2]{RostamiXue18}.}
For nonsymmetric problems, large imaginary spectral components can increase the number of required substeps and reduce the efficiency of these function-specific constructions.

\subsection{The PAIN algorithm}
A generic rational Krylov subspace method with nonorthogonal basis vectors called \emph{poles and interpolation nodes} (PAIN) was proposed in \cite[section 5.4.2]{GuttelThesis} for approximating $f(A)b$. The first basis vector is $v_1 = b/\|b\|_2$, and the subsequent basis {vectors} are $\alpha_j v_{j+1} = (I-A/\xi_j)^{-1}(A-\sigma_j I)v_j$, where the nodes $\sigma_j$ are distinct from the poles $\xi_j$, and the displayed normalization assumes finite nonzero poles outside $\sigma(A)$. Rearranging each recurrence gives $AV_{m+1}\underline{K_m} = V_{m+1}\underline{H_m}$, where $V_{m+1} = [v_1,\ldots,v_{m+1}]$ and 

{\footnotesize
\begin{eqnarray}\label{painKH}
\quad    \underline{K_m} = \left[\begin{array}{cccc}1 & & & \\ {\alpha_1}/{\xi_1} & 1 & & \\ & {\alpha_2}/{\xi_2} & \ddots & \\ & & \ddots & 1 \\ \hline & & & {\alpha_m}/{\xi_m}\end{array}\right], \quad \underline{H_m} = \left[\begin{array}{cccc}\sigma_1 & & & \\ \alpha_1 & \sigma_2 & & \\ & \alpha_2 & \ddots & \\ & & \ddots & \sigma_m \\ \hline & & & \alpha_m\end{array}\right]. 
\end{eqnarray}
}

It can be shown that the approximation $f_m = V_m f(H_m K_m^{-1})\|b\|_2 e_1$ satisfies $f_m = r_m(A)b$, where $r_m(z)$ is a rational function of $(m-1,m-1)$-type that interpolates $f(z)$ at $\{\sigma_j\}_{j=1}^m$ and has poles $\{\xi_j\}_{j=1}^{m-1}$. The main virtue of PAIN is the short-recurrence $f_m = f_{m-1}+c_m v_m$, where $c_m$ is the last element of $f(H_m K_m^{-1})\|b\|_2 e_1$. There is no need to keep the entire basis $V_m$, and there is no basis orthogonalization cost. PAIN requires prescribed nodes and poles, so its practical use depends on useful scalar approximation data. AAA supplies such data adaptively on a prescribed sample set, while selecting that set remains problem-dependent. This requirement also applies to our methods and distinguishes them from adaptive rational Arnoldi \cite{Guttel.Knizhnerman.2013,Gue13b}.

For the same final admissible node and pole sets, Algorithm~\ref{alg:realratlejamv} and PAIN evaluate the same rational interpolant in exact arithmetic. Conjugate pairing can change the pole prefixes used at intermediate steps. The comparisons in section~\ref{sec:experiment} assess the arithmetic realizations of these interpolants.

\subsection{Two-pass Lanczos and truncated Arnoldi methods} For a real symmetric matrix $A$, the standard Lanczos method leads to the relation $AQ_{m}=Q_{m+1}\underline{H_m}$, where $q_1 = b/\|b\|_2$, $Q_{m+1}$ contains an orthonormal basis of the Krylov subspace $\mathcal{K}_{m+1}(A,b)$, and $\underline{H_m}$ is tridiagonal. Similarly, 
the $m+1$-dimensional rational Krylov subspace $\textstyle \mathcal{Q}_{m+1}(A,b) = \mathrm{span}\left\{b,(I-A/\xi_1)^{-1}b,\ldots, \prod_{k=1}^m(I-A/\xi_k)^{-1}b\right\}$ can be generated by a short-recurrence rational Lanczos algorithm. Let $\xi_0 = \infty$ so that $\xi_0^{-1} = 0$. Rational Lanczos generates the relation $AQ_{m+1}\left(\underline{I_m}+D_m \underline{H_m}\right) = Q_{m+1}\underline{H_m}$, where $Q_{m+1}$ contains an orthonormal basis of $\mathcal{Q}_{m+1}$, $\underline{I_m}, \underline{H_m} \in \mathbb{R}^{(m+1) \times m}$ are defined by padding a row below the identity $I_m$ and $H_m$, and $D_m = \mathrm{diag}(\xi_0^{-1},\xi_1^{-1},\ldots,\xi_m^{-1})$. 

In both polynomial and rational Lanczos, let $J_m=Q_m^TAQ_m$ denote the orthogonal compression. The projected approximation is $f_m=Q_mf(J_m)e_1\|b\|_2$. In the polynomial case, $J_m=H_m$; in the rational case, the pencil factors require the appropriate endpoint correction when forming $J_m$, and $J_m$ should not in general be identified with $H_mK_m^{-1}$. The compressed matrix can be obtained from the recurrence data; see \cite[section 2]{PalittaEtal22}. A two-pass implementation forms the small problem in the first pass and then regenerates the basis to accumulate $Q_my_m$, where $y_m=f(J_m)e_1\|b\|_2$. This avoids retaining the full basis, at the cost of a second basis-generation pass.

In a recent study \cite{PalittaEtal22}, the rational Lanczos algorithm has been improved in efficiency and functionality to approximate the bilinear form $u^Tf(A)v$ and the quadratic form $v^Tf(A)v$, without the need to store the entire basis vectors $Q_m$ or double the arithmetic cost; see \cite[sections 2 and 3]{PalittaEtal22} for details.  

For nonsymmetric matrices, sketching can be combined with truncated Arnoldi and a two-pass strategy \cite{GuttelSchweitzer23,PalittaSchweitzerSimoncini25}. Only a window of recent large vectors is retained during basis generation; the reduced problem uses sketches, and the final action is accumulated during a second pass. The orthogonalization window controls the local orthogonalization work, but a short window need not produce a well-conditioned full basis. Section~\ref{subsec:polytest} reports a comparison in which a successful window was selected separately for each matrix and function. This provides a useful contrast with interpolation: the latter uses prescribed scalar approximation data instead of selecting an orthogonalization window, and accumulates the action in one pass.

{\bf Summary.} These methods make different choices about approximation information, orthogonalization, and basis regeneration. Our focus is on the arithmetic realization of prescribed interpolation: preserving conjugation makes real matrix operations possible without retaining a full basis.

\section{Why preserve conjugation symmetry?}\label{sec:asymfac}
The main reason for preserving conjugation is algebraic: conjugation-preserving data admit real interpolation coefficients in a suitable basis, allowing the matrix-vector stage to use real arithmetic. Potential theory gives a complementary observation. When both plates of a condenser are symmetric with respect to the real axis, their equilibrium distribution has the same symmetry. This observation makes conjugate pairing a natural design choice, but does not make exact pairing necessary for every asymptotically good sequence.

For this discussion, let $\Sigma$ and $\Xi$ be disjoint compact sets of positive logarithmic capacity, with positive distance between them and with $\overline{\Sigma}=\Sigma$ and $\overline{\Xi}=\Xi$. Let $\mathcal M(\Sigma)$ and $\mathcal M(\Xi)$ denote the probability measures supported on the respective sets. The condenser equilibrium measure $\mu^*=\mu^*_{\Sigma}-\mu^*_{\Xi}$ is the unique finite-energy minimizer of
\[
 I(\mu)=\iint\log\frac{1}{|z-\zeta|}\,d\mu(z)\,d\mu(\zeta),
 \qquad \mu\in\mathcal M(\Sigma)-\mathcal M(\Xi).
\]
Writing $V=I(\mu^*)>0$, the condenser capacity is $\operatorname{cap}(\Sigma,\Xi)=1/V$; see \cite[Chapter VIII]{SaffTotik1997} and \cite[section 7.2]{GuttelThesis}.

\begin{lemma}\label{lem:equilibrium_conjugation}
Let $C(z)=\overline z$, and define the reflected measure by $(C_\#\mu)(B)=\mu(C^{-1}(B))$ for every Borel set $B$. Then $C_\#\mu^*=\mu^*$.
\end{lemma}
\begin{proof}
Reflection preserves each plate, so $C_\#\mu^*$ is admissible. Since $|C(z)-C(\zeta)|=|z-\zeta|$, a change of variables gives $I(C_\#\mu^*)=I(\mu^*)$. Uniqueness of the minimizer therefore implies $C_\#\mu^*=\mu^*$.
\end{proof}

Suppose that normalized signed node-and-pole counting measures $\mu_m$ converge weak-star to $\mu^*$. For every fixed continuous test function $g$ on $\Sigma\cup\Xi$,
\[
 \int g\,d(\mu_m-C_\#\mu_m)
 =\int (g-g\circ C)\,d\mu_m
 \longrightarrow \int (g-g\circ C)\,d\mu^*=0.
\]
Thus their limiting distributions are symmetric. The finite sets themselves need not be exactly invariant under conjugation. Our constructions enforce this stronger finite-degree property to obtain real coefficients, as proved in sections~\ref{sec:poly} and \ref{sec:lmrksm}.

Classical potential theory gives capacity-based convergence estimates for regular condenser geometries and suitably distributed nodes and poles \cite[section~7.2]{GuttelThesis}. Such estimates require additional analyticity and distribution hypotheses. We use symmetry here as a design principle, not as a claim that conjugate-completed Leja or AAA data attain an optimal asymptotic rate. The finite-degree interpolation results below are independent of that question.

\section{Limited-memory polynomial interpolation}\label{sec:poly}
We first construct a polynomial interpolant $p_m$ of degree at most $m$ on a conjugation-invariant region $D$. The key is to treat a real interpolation point as a one-point block and a nonreal conjugate pair as a two-point block. This gives a Newton-like basis with real coefficients, while preserving the progressive structure of interpolation.

\begin{definition}\label{def:conjordered}
An ordered set of distinct points is \emph{conjugation-invariant} if each point is real or belongs to a pair of adjacent nonreal conjugate points. Thus the ordering consists of real one-point blocks and conjugate two-point blocks.
\end{definition}

Set $p_{-1}=0$ and $q_0=1$. Suppose that the completed blocks contain $z_0,\ldots,z_{j-1}$, that $p_{j-1}$ interpolates $f$ there, and that $q_j$ vanishes at these points. Choose a new point by the Leja criterion
\[
 z_j\in\mathop{\rm argmax}_{z\in\partial D}\prod_{\ell=0}^{j-1}|z-z_\ell|.
\]
In practice, the maximization is performed on a finite conjugation-invariant sample set $Z$ of the boundary. If $z_j$ is real, append it and choose $c_j$ so that $p_j=p_{j-1}+c_jq_j$ interpolates at $z_j$ as well. If $z_j$ is nonreal, append $z_j$ and $z_{j+1}=\overline{z_j}$ together. Introduce a real number $\sigma_j$ and set
\[
 q_{j+1}(z)=\nu_{j+1}(z-\sigma_j)q_j(z),\qquad
 p_{j+1}=p_{j-1}+c_jq_j+c_{j+1}q_{j+1}.
\]
Both new basis functions vanish at the preceding nodes. The two coefficients are chosen to interpolate at the new pair, so all preceding interpolation conditions remain satisfied. Appending the conjugate point is a symmetry-completion step; that second point need not itself maximize the next ordinary Leja criterion.

The complete basis is defined by $q_0=1$ and, for $j\geq1$,
\begin{equation}\label{defpolybasisnew}
\begin{aligned}
 q_j(z)&=\begin{cases}
 \displaystyle\gamma_j(z-\sigma_{j-1})\prod_{\ell=0}^{j-2}(z-z_\ell),
 &z_{j-1}\notin\mathbb R,\ z_j=\overline{z_{j-1}},\\[5pt]
 \displaystyle\gamma_j\prod_{\ell=0}^{j-1}(z-z_\ell),&\text{otherwise},
 \end{cases}\\[4pt]
 &=\begin{cases}
 \nu_j(z-\sigma_{j-1})q_{j-1}(z),
 &z_{j-1}\notin\mathbb R,\ z_j=\overline{z_{j-1}},\\
 \nu_j(z-z_{j-1})(z-z_{j-2})q_{j-2}(z),
 &z_{j-1}=\overline{z_{j-2}}\notin\mathbb R,\\
 \nu_j(z-z_{j-1})q_{j-1}(z),&z_{j-1}\in\mathbb R.
 \end{cases}
\end{aligned}
\end{equation}
Here $\gamma_0=1$ and the nonzero real $\gamma_j$ are cumulative normalizations. The local recurrence factors are $\nu_j=\gamma_j/\gamma_{j-1}$ in the first and third branches and $\nu_j=\gamma_j/\gamma_{j-2}$ in the second. We choose them so that $\|q_j(Z)\|_\infty=1$. This controls the sampled scale of the basis functions and reduces the risk of overflow, underflow, and poorly scaled interpolation columns; it does not by itself control the condition number of the interpolation system.

Every $q_j$ has real coefficients: its factors consist of real points, complete conjugate pairs, and possibly $z-\sigma_{j-1}$. The choice of $\sigma_j\in\mathbb R$ does not affect the final interpolant in exact arithmetic. In the polynomial experiments we use $\sigma_j=\operatorname{mean}(\operatorname{Re}Z)$. Tests with $\sigma_j=0$, $\sigma_j=\mu\operatorname{Re}(z_j)$ for $0.1\leq\mu\leq10$, and $\sigma_j=\|A\|_\infty$ did not give a consistent improvement in accuracy.

The coefficients must use the same normalization as the basis. In particular, the paired update is $c_jq_j+c_{j+1}q_{j+1}$, with $q_{j+1}=\nu_{j+1}(z-\sigma_j)q_j$. An equivalent implementation can absorb $\nu_{j+1}$ into the coefficient of the unnormalized second column. We initialize the algorithms with a real boundary point, as available in the regions used here; a conjugate-pair initialization is also possible. The final index $m$ always completes a block. A provisional degree budget is increased by one if necessary to avoid splitting a pair, and no further basis update is needed after the final interpolant is formed.

As an example, suppose that $z_0$ and $z_1 = \overline{z_0}$ are complex, $z_2$ is real, $z_3$ and $z_4 = \overline{z_3}$ are complex, and $z_5$ and $z_6 = \overline{z_5}$ are also complex. To construct $p_6(z)$ that interpolates $f(z)$ at $\{z_0,z_1,\ldots,z_6\}$, we obtain the coefficients $c_0,\ldots,c_6$ by solving the following quasi-lower triangular linear system (for this illustration, set all $\sigma_j=0$ and $\nu_1=\nu_4=\nu_6=1$)

{\footnotesize
\begin{eqnarray}\label{polylinex}
   \qquad  \left[\begin{array}{ccccccc} q_0(z_0) & z_0 q_0(z_0) & & & & & \\ q_0(\overline{z_0}) & \overline{z_0} q_0(\overline{z_0}) & & & & & \\ q_0(z_2) & z_2 q_0(z_2) & q_2(z_2) & & & & \\ q_0(z_3) & z_3 q_0(z_3) & q_2(z_3) & q_3(z_3) & z_3 q_3(z_3) & & \\ q_0(\overline{z_3}) & \overline{z_3} q_0(\overline{z_3}) & q_2(\overline{z_3}) & q_3(\overline{z_3}) & \overline{z_3} q_3(\overline{z_3}) & & \\ q_0(z_5) & z_5 q_0(z_5) & q_2(z_5) & q_3(z_5) & z_5 q_3(z_5) & q_5(z_5) & z_5 q_5(z_5) \\ q_0(\overline{z_5}) & \overline{z_5} q_0(\overline{z_5}) & q_2(\overline{z_5}) & q_3(\overline{z_5}) & \overline{z_5} q_3(\overline{z_5}) & q_5(\overline{z_5}) & \overline{z_5} q_5(\overline{z_5}) \end{array} \right]\left[\begin{array}{c} c_0 \\ c_1 \\ c_2 \\ c_3 \\ c_4 \\ c_5 \\ c_6 \end{array}\right] = \left[\begin{array}{c} f(z_0) \\ f(\overline{z_0}) \\ f(z_2) \\ f(z_3) \\ f(\overline{z_3}) \\ f(z_5) \\ f(\overline{z_5})\end{array}\right],
\end{eqnarray}
}
where the basis polynomials are
\begin{eqnarray}
    && q_0(z) = 1, \quad q_1(z) = \nu_1 zq_0(z), \quad q_2(z) = \gamma_2(z-z_0)(z-\overline{z_0}), \\ \nonumber
    && q_3(z) = \gamma_3 (z-z_0)(z-\overline{z_0})(z-z_2),\quad q_4(z) = \nu_4 zq_3(z) \\ \nonumber
    && q_5(z) = \gamma_5(z-z_0)(z-\overline{z_0})(z-z_2)(z-z_3)(z-\overline{z_3}), \quad q_6(z) = \nu_6 zq_5(z).
\end{eqnarray}

\begin{algorithm}
	\caption{Scalar polynomial interpolation at conjugate-completed Leja points}
	\label{alg:realpolylejas}
	\begin{algorithmic}[1]
	    \Require degree budget $m$; a conjugation-invariant grid $Z\subset\partial D$ with sufficiently many points and a real point; $f(\overline z)=\overline{f(z)}$ on $D=\overline D$.
        \Ensure a polynomial approximation $p_m(z) \approx f(z)$ on $D$;
	    \State let $j=0$; define $p_{-1}(z) = 0$, and $q_0(z)\equiv 1$;
        \While{$j \leq m$}
	    \State choose an unused $z_j\in Z$ maximizing $\prod_{i=0}^{j-1}|z-z_i|$; if $j=0$, choose a real $z_0\in Z$; 
        \If{$z_j$ is real}
            \State find $c_j \in \mathbb{R}$ s.t. $p_j(z) = p_{j-1}(z) + c_j q_j(z)$ interpolates $f(z)$ at $z_0,\ldots,z_j$;
            \State $j = j+1$;
        \Else
            \State if $j=m$, set $m\leftarrow m+1$; let $z_{j+1} = \overline{z_j}$, choose $\sigma_j\in\mathbb{R}$, and form $q_{j+1}(z)=\nu_{j+1}(z-\sigma_j)q_j(z)$ with the chosen normalization;
            \State find $c_j,c_{j+1}\in\mathbb{R}$ such that $p_{j+1}=p_{j-1}+c_jq_j+c_{j+1}q_{j+1}$ interpolates $f$ at $z_0,\ldots,z_{j+1}$;
            \State $j=j+2$;
        \EndIf
        \State if $j\le m$, update $q_j(z)$ by \eqref{defpolybasisnew}, with $\|q_j(Z)\|_\infty=1$; 
        \EndWhile
    \end{algorithmic}
\end{algorithm}

The system is solved progressively by block forward substitution. For the example above, one first solves for $c_0,c_1$, then for $c_2$, and then for $c_3,c_4$ and $c_5,c_6$. The implementation applies local iterative refinement to each new block. At a completed real node the update is $p_j=p_{j-1}+c_jq_j$; at the end of a conjugate pair it is $p_j=p_{j-2}+c_{j-1}q_{j-1}+c_jq_j$. The following result explains why these updates require real coefficients only.

\begin{theorem}\label{thmcoefspolylinex}
    Consider the construction of $p_m(z)$ with conjugation-invariant distinct interpolation points $z_0,z_1,\ldots,z_m$ by Algorithm \ref{alg:realpolylejas}. The quasi-lower triangular linear system exemplified in \eqref{polylinex} has a unique real solution. 
\end{theorem}
\begin{proof}
The coefficient matrix is block lower triangular. At a real node, its diagonal entry is $q_j(z_j)\ne0$, because the nodes are distinct. At a conjugate pair $z_j,z_{j+1}=\overline{z_j}$, the diagonal block is
\[
 \begin{bmatrix}
 q_j(z_j)&\nu_{j+1}(z_j-\sigma_j)q_j(z_j)\\
 \overline{q_j(z_j)}&\nu_{j+1}(\overline{z_j}-\sigma_j)\overline{q_j(z_j)}
 \end{bmatrix}.
\]
Its determinant is $\nu_{j+1}(\overline{z_j}-z_j)|q_j(z_j)|^2\ne0$. Thus the full matrix $Q$ is nonsingular. Let $J$ exchange the rows associated with conjugate nodes and fix those associated with real nodes. Since the basis polynomials have real coefficients, $\overline Q=JQ$; the data vector $f_Z$ satisfies $\overline{f_Z}=Jf_Z$. Conjugating $Qc=f_Z$ therefore gives $Q\overline c=f_Z$. Uniqueness implies $c=\overline c$.
\end{proof}

The real coefficients can also be read directly from a single complex residual. At the start of a nonreal block, put
\[
 \eta_j=\frac{f(z_j)-p_{j-1}(z_j)}{q_j(z_j)}.
\]
The equation $c_j+\nu_{j+1}(z_j-\sigma_j)c_{j+1}=\eta_j$ gives
\begin{equation}\label{eq:pairedcoefficients}
 c_{j+1}=\frac{\operatorname{Im}\eta_j}{\nu_{j+1}\operatorname{Im}z_j},\qquad
 c_j=\operatorname{Re}\eta_j-
       \frac{\operatorname{Re}z_j-\sigma_j}{\operatorname{Im}z_j}
       \operatorname{Im}\eta_j.
\end{equation}
This is an exact-arithmetic description of the two-point block; it also makes clear that nearly real pairs can require care in floating-point arithmetic.

\begin{corollary}\label{col:polybasisindp}
    The polynomials $\{q_j(z)\}_{j=0}^m$ defined in \eqref{defpolybasisnew} form a basis of $\mathbb{P}_m$, the set of all polynomials of degree $\leq m$.
\end{corollary}
\begin{proof}
Each $q_j$ has degree exactly $j$ and a nonzero leading coefficient. Polynomials of distinct successive degrees are linearly independent, so the $m+1$ polynomials form a basis of $\mathbb P_m$.
\end{proof}

\begin{algorithm}
	\caption{Polynomial matrix-function action at conjugate-completed Leja points}
	\label{alg:realpolylejamv}
	\begin{algorithmic}[1]
	    \Require the scalar inputs of Algorithm~\ref{alg:realpolylejas}; $A\in\mathbb R^{n\times n}$ and $b\in\mathbb R^n$, with $W(A)\subset D$.
        \Ensure a polynomial $p_m(z) \approx f(z)$ on $D$, and the vector $p_m(A)b \approx f(A)b$;
	    \State let $j=0$; define $p_{-1}(z) = 0$ (scalar), $p_{-1}(A)b = 0$ (vector), $q_0(z)\equiv 1$ and $q_0(A)b \equiv b$;
        \While{$j \leq m$}
	    \State choose an unused $z_j\in Z$ maximizing $\prod_{i=0}^{j-1}|z-z_i|$; if $j=0$, choose a real $z_0\in Z$; 
        \If{$z_j$ is real}
            \State find $c_j \in \mathbb{R}$ s.t. $p_j(z) = p_{j-1}(z) + c_j q_j(z)$ interpolates $f(z)$ at $z_0,\ldots,z_j$;
            \State evaluate $p_j(A)b = p_{j-1}(A)b + c_j q_j(A)b$;
            \State $j = j+1$;
        \Else
            \State if $j=m$, set $m\leftarrow m+1$; let $z_{j+1} = \overline{z_j}$, choose $\sigma_j\in\mathbb{R}$, and form $q_{j+1}(z)=\nu_{j+1}(z-\sigma_j)q_j(z)$ with the chosen normalization;
            \State find $c_j,c_{j+1}\in\mathbb{R}$ such that $p_{j+1}=p_{j-1}+c_jq_j+c_{j+1}q_{j+1}$ interpolates $f$ at $z_0,\ldots,z_{j+1}$;
            \State evaluate $p_{j+1}(A)b=p_{j-1}(A)b+c_jq_j(A)b+c_{j+1}\nu_{j+1}(A-\sigma_jI)q_j(A)b$;
            \State $j=j+2$;
        \EndIf
        \State if $j\le m$, update $q_j(z)$ by \eqref{defpolybasisnew}, with $\|q_j(Z)\|_\infty=1$; 
        \State if $j\le m$, update $q_j(A)b$ by the same recurrence in real arithmetic;
        \EndWhile
    \end{algorithmic} 
\end{algorithm} 

\paragraph{Evaluation of the matrix action.}
Once the scalar interpolation data are available, all large-vector operations in Algorithm~\ref{alg:realpolylejamv} use real arithmetic. A real node gives
\[
 q_j(A)b=\nu_j(A-z_{j-1}I)q_{j-1}(A)b.
\]
At the end of a conjugate block,
\[
 q_j(A)b=\nu_j\bigl(A^2-2\operatorname{Re}(z_{j-1})A+|z_{j-1}|^2I\bigr)q_{j-2}(A)b.
\]
The expression is evaluated by successive real matrix--vector products, without forming $A^2$. The intermediate column is obtained from $\nu_{j+1}(A-\sigma_jI)q_j(A)b$. Only the current block, the required predecessor, and the accumulated approximation need to be retained.

The difference from real Leja interpolation in section~\ref{subsec:realleja} is the location of the nodes. That method maps precomputed nodes to the focal interval of an ellipse, whereas Algorithm~\ref{alg:realpolylejamv} selects conjugate-completed nodes on a sampled boundary of $D$. These choices need not have the same convergence rate. In the exponential examples of section~\ref{subsec:polytest}, the boundary construction sometimes uses about twice as many interpolation points as the focal-interval construction.

\section{Limited-memory shift-invert interpolation}\label{sec:lmsira} Polynomial approximation can be slow when the region is large or lies close to a singularity of $f$. A fixed repeated pole can improve the scalar approximation while allowing a single shifted factorization to be reused. For the exponential action $\exp(-hA)b$, a common degree-dependent choice is $s=-m/h$, motivated by shift-invert convergence bounds \cite{Moret.Novati.2004,Borwein1983}. Repeated-pole choices for Cauchy--Stieltjes functions are discussed in \cite[section~6.1]{Beckermann2009}. Here we transfer the polynomial construction to a fixed real pole.

We choose a real shift $s \not\in D$, and define the transformed domain $E = \{w\,:\,w = \frac{1}{z-s}, z \in D\}$. Then we define $g(w) := f\left(w^{-1}+s\right)=f(z) $. If $p_m(w)$ interpolates $g$ on $E$, then $p_m\left((z-s)^{-1}\right) = p_m(w) \approx g(w) = f(z)$ is a shift-invert (single-pole) rational approximation to $f(z)$. Let $w_0,\ldots,w_m \in E$ be the conjugation-invariant Leja interpolation points found in the construction of $p_m(w) \approx g(w)$. The shift-invert rational approximation to $f(A)b$ is outlined in Algorithm \ref{alg:realpolylejasimv}, obtained by replacing $z$ with $w$, $f(z)$ with $g(w)$, and matrix $A$ with $(A-sI)^{-1}$ in Algorithm~\ref{alg:realpolylejamv}. The real-node and conjugate-pair updates again give short recurrences. 

The substitution is an identity of rational functional calculus. It does not
require $W((A-sI)^{-1})\subset E$, which need not hold for a nonsymmetric matrix.
Indeed, $q(z)=p_m((z-s)^{-1})$ has no poles on $D$, and
\eqref{crouzeixbound} applies directly to $f-q$ on the original set containing
$W(A)$. Thus a uniform bound for $g-p_m$ on $E$ provides the required scalar
bound for $f-q$ on $D$.

\begin{algorithm}
	\caption{Shift-invert matrix-function action at conjugate-completed Leja points}
	\label{alg:realpolylejasimv}
	\begin{algorithmic}[1]
	    \Require the scalar inputs of Algorithm~\ref{alg:realpolylejas}; $A\in\mathbb R^{n\times n}$ and $b\in\mathbb R^n$, with $W(A)\subset D$.
        \Ensure $p_m((z-s)^{-1}) \approx f(z)$ on $D$, and the vector $p_m(M)b \approx f(A)b$;
        \State choose $s\in\mathbb R\setminus D$; set $E=\{(z-s)^{-1}:z\in D\}$, $Z_E=\{(z-s)^{-1}:z\in Z\}$, and $g(w)=f(w^{-1}+s)$;
	    \State let $M=(A-sI)^{-1}$ denote a solve operator (do not form the inverse); set $j=0$; define $p_{-1}(w) = 0$ (scalar), $p_{-1}(M)b = 0$ (vector), $q_0(w)\equiv 1$ and $q_0(M)b \equiv b$;
        \While{$j \leq m$}
	    \State choose an unused $w_j\in Z_E$ maximizing $\prod_{i=0}^{j-1}|w-w_i|$; if $j=0$, choose a real $w_0\in Z_E$; 
        \If{$w_j$ is real}
            \State find $c_j \in \mathbb{R}$ s.t. $p_j(w) = p_{j-1}(w) + c_j q_j(w)$ interpolates $g(w)$ at $w_0,\ldots,w_j$;
            \State evaluate $p_j(M)b = p_{j-1}(M)b + c_j q_j(M)b$;
            \State $j = j+1$;
        \Else
            \State if $j=m$, set $m\leftarrow m+1$; let $w_{j+1}=\overline{w_j}$, choose $\sigma_j\in\mathbb{R}$, and form $q_{j+1}(w)=\nu_{j+1}(w-\sigma_j)q_j(w)$ with the chosen normalization;
            \State find $c_j,c_{j+1}\in\mathbb{R}$ such that $p_{j+1}=p_{j-1}+c_jq_j+c_{j+1}q_{j+1}$ interpolates $g$ at $w_0,\ldots,w_{j+1}$;
            \State evaluate $p_{j+1}(M)b = p_{j-1}(M)b+c_jq_j(M)b$\\ \qquad\qquad $+c_{j+1}\nu_{j+1}(M-\sigma_j I)q_j(M)b$;
            \State $j=j+2$;
        \EndIf
        \State if $j\le m$, update $q_j(w)$ by \eqref{defpolybasisnew}, normalized on $Z_E$; 
        \State if $j\le m$, update $q_j(M)b$ by the same recurrence;
        \EndWhile
    \end{algorithmic} 
\end{algorithm}

\section{Limited-memory multi-pole interpolation}\label{sec:lmrksm}
In this section, we construct a rational function $r_m$ of type at most $(m,m)$
that interpolates $f$ at distinct points $Z_{m+1}=\{z_0,\ldots,z_m\}\subset D$,
with prescribed poles $\Xi_m=\{\xi_0,\ldots,\xi_{m-1}\}$.
The node and pole sets are conjugation-invariant ordered sets, and $z_0$ is real.
We assume that the poles lie outside $D$, so that they do not coincide with
interpolation points or eigenvalues of $A$. These conditions give
$r_m(\overline z)=\overline{r_m(z)}$ and permit evaluation of $r_m(A)b$.

\paragraph{Selecting conjugation-preserving data.}
To select the nodes and poles, we use a paired-support variant of AAA
\cite{NakatsukasaEtal18}. Let $Z$ be a finite conjugation-invariant sample set
representing $D$, with at least one real point. Initially, choose
$z_0\in\mathop{\rm argmax}_{z\in Z\cap\mathbb R}|f(z)|$ and set
$r_0(z)=f(z_0)$. At step $m\ge1$, choose an unused point maximizing
$|f(z)-r_{m-1}(z)|$. If this point $z_m$ is real, append it and put $p=0$;
otherwise, append both $z_m$ and $z_{m+1}=\overline{z_m}$ and put $p=1$.
For the resulting support points $z_0,\ldots,z_{m+p}$, define
\[
 r_{m+p}(z)=
 \frac{\displaystyle\sum_{j=0}^{m+p}\frac{w_jf(z_j)}{z-z_j}}
 {\displaystyle\sum_{j=0}^{m+p}\frac{w_j}{z-z_j}}.
\]
Writing the remaining sample points as $\zeta_1,\ldots,\zeta_L$, AAA chooses
a unit right singular vector associated with the smallest singular value of
the Loewner matrix
\[
 (\mathcal L_m)_{ij}=
 \frac{f(\zeta_i)-f(z_j)}{\zeta_i-z_j},
 \qquad 1\le i\le L,\quad 0\le j\le m+p.
\]
This minimizes the Euclidean norm of the linearized residual at the unused
sample points. It does not by itself certify a uniform error on the continuous
domain $D$ or the exclusion of poles from $D$.
For nonzero support weights, the barycentric formula interpolates the data
at its support points. Its denominator zeros are obtained from the pencil
\begin{eqnarray}\label{gevpole}
 \left(
 \left[\begin{array}{ccccc}
 0&w_0&w_1&\ldots&w_{m+p}\\
 1&z_0&&&\\1&&z_1&&\\\vdots&&&\ddots&\\1&&&&z_{m+p}
 \end{array}\right],
 \left[\begin{array}{ccccc}
 0&&&&\\&1&&&\\&&1&&\\&&&\ddots&\\&&&&1
 \end{array}\right]\right).
\end{eqnarray}
Generically, this pencil has two infinite eigenvalues and $m+p$ finite
denominator zeros. After each accepted support block, the support index is increased accordingly and the greedy process continues until its sampled stopping criterion is met. The fixed-pole construction below covers the nondegenerate case with $m+1$ distinct support points, nonzero weights, and $m$ distinct finite poles outside $D$. Degree reduction, repeated poles, or exact cancellations require a separate treatment; AAA alone does not guarantee these nondegeneracy and admissibility conditions.

Pairing the support points allows a symmetric minimizing weight vector to be chosen. An arbitrary singular vector need not have that symmetry: its phase is undetermined, and a multiple smallest singular value also leaves a choice within the minimizing subspace. The following statement makes the distinction precise.

\begin{proposition}\label{prop:aaa_symmetry}
Suppose that the support and remaining sample sets are invariant under
conjugation and $f(\overline z)=\overline{f(z)}$ on these sets. There exists a
unit vector minimizing $\|\mathcal L_mw\|_2$ whose weights satisfy
$w_j\in\mathbb R$ at real support points and
$w_k=\overline{w_j}$ when $z_k=\overline{z_j}$.
If the smallest right singular subspace is one-dimensional, every minimizing
vector has this form after multiplication by a constant of modulus one.
For such a choice, the barycentric rational function is conjugation-preserving,
and its poles, after cancellations, form a conjugation-invariant multiset.
\end{proposition}
\begin{proof}
Let $J_r$ and $J_c$ be the permutation matrices that exchange conjugate
remaining sample points and conjugate support points, respectively, while
fixing the real points. The divided differences satisfy
$\mathcal L_mJ_c=J_r\overline{\mathcal L_m}$.
For each conjugate pair use the unitary block
\[
 U_2=\frac1{\sqrt2}\begin{bmatrix}1&i\\1&-i\end{bmatrix},
\]
and use the scalar $1$ for each real point. The resulting column and row
transformations $U_c,U_r$ satisfy $J_c\overline{U_c}=U_c$ and
$J_r\overline{U_r}=U_r$. Therefore
$\widehat{\mathcal L}_m=U_r^*\mathcal L_mU_c$ is real. It has the same
singular values as $\mathcal L_m$ and admits a real unit minimizing vector
$v$. Then $w=U_cv$ is a minimizing vector with the stated conjugacy.
If the minimum right singular subspace is one-dimensional, all its unit
vectors differ by a phase, proving the second assertion.
Finally, conjugating and permuting the terms in the barycentric numerator
and denominator preserves both sums in the required way. Multiplication by
$\prod_j(z-z_j)$ gives real-coefficient numerator and denominator polynomials,
so the uncancelled poles occur in conjugate pairs.
\end{proof}

The transformed real SVD gives a structure-preserving way to select the weights, including when the smallest singular value is multiple. This exact-arithmetic result concerns the weights, not the conditioning of pole extraction. The experiments compare matrix-action realizations; they do not time or compare alternative scalar SVD implementations.

\paragraph{A real short-recurrence basis.}

We next describe the limited-memory evaluation for fixed admissible nodes and
poles. It is useful to count poles rather than use the last pole index.
Let $d_0=0$. For $1\le j\le m$, define $d_j$ as follows. If
$z_{j-1}\notin\mathbb R$ and $z_j=\overline{z_{j-1}}$, set
$d_j=d_{j-1}$. In every other case, let $d_j$ be the smallest integer
$d\in\{0,\ldots,m\}$ with
\[
 d\ge\min(j+1,m),\qquad
 \{\xi_0,\ldots,\xi_{d-1}\}\ \hbox{conjugation-invariant}.
\]
This integer always exists because the full pole set is conjugation-invariant.
In particular, $d_j\ge j$, $d_m=m$, and the pole prefix is unchanged between
the two columns associated with a conjugate node pair.
Set $P_j(z)=\prod_{\ell=0}^{j-1}(z-z_\ell)$ and
$\ell_d(z)=\prod_{\ell=0}^{d-1}(z-\xi_\ell)$, with empty products equal to one.
With nonzero real normalization constants $\gamma_j$ and $\gamma_0=1$, define
\begin{equation}\label{defratbasis}
 q_0(z)=1,\qquad
 q_j(z)=\begin{cases}
 \displaystyle\gamma_j\frac{(z-\sigma_{j-1})P_{j-1}(z)}{\ell_{d_j}(z)},
 &z_{j-1}\notin\mathbb R,\ z_j=\overline{z_{j-1}},\\[6pt]
 \displaystyle\gamma_j\frac{P_j(z)}{\ell_{d_j}(z)},&\text{otherwise},
 \end{cases}
\end{equation}
where every $\sigma_j$ is real. The paired case has
$q_j=(\gamma_j/\gamma_{j-1})(z-\sigma_{j-1})q_{j-1}$.
Thus every $q_j$ is conjugation-preserving and has numerator degree $j$ and
denominator degree $d_j$. Its common-denominator numerator with denominator
$\ell_m$ has degree at most $m$.

Let $\nu_j$ denote the ratio of cumulative normalizations associated with
the preceding basis function in each branch below. The short recurrences are
\begin{equation}\label{defratbasisnew}
 q_j(z)=\begin{cases}
 \nu_j(z-\sigma_{j-1})q_{j-1}(z),
 &z_{j-1}\notin\mathbb R,\ z_j=\overline{z_{j-1}},\\[5pt]
 \displaystyle\nu_j
 \frac{(z-z_{j-1})(z-z_{j-2})}
 {\prod_{\ell=d_{j-2}}^{d_j-1}(z-\xi_\ell)}q_{j-2}(z),
 &z_{j-1}=\overline{z_{j-2}}\notin\mathbb R,\\[9pt]
 \displaystyle\nu_j\frac{z-z_{j-1}}
 {\prod_{\ell=d_{j-1}}^{d_j-1}(z-\xi_\ell)}q_{j-1}(z),
 &z_{j-1}\in\mathbb R.
 \end{cases}
\end{equation}
Here $\nu_j=\gamma_j/\gamma_{j-1}$ in the first and third branches, and
$\nu_j=\gamma_j/\gamma_{j-2}$ in the second. Every denominator in the recurrence
contains only newly added poles; an empty product is one. This convention
also handles the last basis function after all poles have been used.
The local normalizations are selected so that
$\|q_j\|_D\approx\|q_j(Z)\|_\infty=1$.
We take $\sigma_{j-1}=\operatorname{Re}(z_{j-1})$ by default. As discussed in
section~\ref{sec:experiment}, the conditioning of the interpolation system
also depends on the ordering of nodes and poles.

For illustration, let $z_0,z_1,z_6$ be real, and let
$z_3=\overline{z_2}$ and $z_5=\overline{z_4}$ be two nonreal conjugate pairs.
Let $\xi_1=\overline{\xi_0}$ and $\xi_4=\overline{\xi_3}$ be nonreal pairs,
with $\xi_2,\xi_5$ real. Then $r_6$ is formed from
\begin{align*}
 q_0(z)&=1,&q_1(z)&=\gamma_1\frac{z-z_0}{\ell_2(z)},\\
 q_2(z)&=\gamma_2\frac{(z-z_0)(z-z_1)}{\ell_3(z)},&
 q_3(z)&=\nu_3(z-\sigma_2)q_2(z),\\
 q_4(z)&=\gamma_4\frac{(z-z_0)(z-z_1)(z-z_2)(z-\overline{z_2})}{\ell_5(z)},&
 q_5(z)&=\nu_5(z-\sigma_4)q_4(z),\\
 q_6(z)&=\gamma_6\frac{\prod_{\ell=0}^5(z-z_\ell)}{\ell_6(z)}.
 \end{align*}
The interpolation coefficients solve the quasi-lower triangular system
\begin{equation}\label{ratlinex}
 Qc=f_Z,\qquad
 Q_{ij}=q_j(z_i),\quad(f_Z)_i=f(z_i),\qquad 0\le i,j\le m.
\end{equation}
Its diagonal blocks have size one for a real node and size two for a
conjugate pair. The normalization factors and real shifts are included in $Q$. The coefficients can be computed progressively by one- or two-point block solves; formula~\eqref{eq:pairedcoefficients} applies to a paired block with the current rational basis and residual.

\begin{algorithm}
 \caption{Limited-memory rational Krylov evaluation with prescribed nodes and poles}
 \label{alg:realratlejamv}
 \begin{algorithmic}[1]
 \Require $f(\overline z)=\overline{f(z)}$ on $D=\overline D$; conjugation-invariant ordered distinct nodes $z_0,\ldots,z_m\in D$ with $z_0$ real, and poles $\xi_0,\ldots,\xi_{m-1}\notin D$; $A\in\mathbb R^{n\times n}$ with $W(A)\subset D$, and $b\in\mathbb R^n\setminus\{0\}$.
 \Ensure the interpolant $r_m(z_i)=f(z_i)$ and its action $r_m(A)b$.
 \State form the scalar bases \eqref{defratbasisnew} and solve \eqref{ratlinex} for $c_0,\ldots,c_m$; set $j=0$, $r_{-1}(A)b=0$, $q_0(A)b=b$;
 \While{$j\le m$}
 \If{$z_j$ is real}
 \State $r_j(A)b=r_{j-1}(A)b+c_jq_j(A)b$;
 \State $j\leftarrow j+1$;
 \Else
 \State $q_{j+1}(A)b=\nu_{j+1}(A-\sigma_jI)q_j(A)b$;
 \State $r_{j+1}(A)b=r_{j-1}(A)b+c_jq_j(A)b+c_{j+1}q_{j+1}(A)b$;
 \State $j\leftarrow j+2$;
 \EndIf
 \If{$j\le m$}
 \State obtain $q_j(A)b$ from $q_{j-1}(A)b$ or $q_{j-2}(A)b$ by \eqref{defratbasisnew}, using each newly added real pole or conjugate pair;
 \EndIf
 \EndWhile
 \end{algorithmic}
\end{algorithm}

\begin{theorem}\label{thm:coefsratlinex}
Under the node and pole assumptions of Algorithm~\ref{alg:realratlejamv},
the interpolation system \eqref{ratlinex} has a unique real solution.
\end{theorem}
\begin{proof}
At a real node $z_j$, the diagonal entry is
$q_j(z_j)=\gamma_jP_j(z_j)/\ell_{d_j}(z_j)\ne0$.
If $z_j,z_{j+1}=\overline{z_j}$ form a nonreal pair, then
$q_j(z_j)\ne0$ and $q_{j+1}=\nu_{j+1}(z-\sigma_j)q_j$.
The determinant of the corresponding diagonal block is
\[
 \nu_{j+1}(\overline{z_j}-z_j)|q_j(z_j)|^2\ne0.
\]
All these entries are finite because nodes and poles are disjoint. Hence
$Q$ is nonsingular. Let $J$ exchange conjugate node rows and fix real rows.
The real basis functions satisfy $\overline Q=JQ$, and the data satisfy
$\overline{f_Z}=Jf_Z$. Conjugating $Qc=f_Z$ therefore gives
$Q\overline c=f_Z$. Uniqueness implies $c=\overline c$.
\end{proof}

We solve the scalar interpolation system \eqref{ratlinex} with MATLAB's backslash, apply iterative refinement \cite[section~12.2]{Higham2002}, and take the real part of the result. For AAA tolerances between $10^{-12}$ and $10^{-9}$, the experiments typically produced interpolation-system condition numbers of order $10^4$--$10^7$ after ordering the nodes and poles. The attainable accuracy depends on that conditioning as well as the scalar approximation error. At tighter tolerances, inaccurate poles and numerical Froissart doublets can become significant \cite{Froissart1969,Beckermann2018,NakatsukasaEtal18}; conjugation symmetry alone does not prevent these effects.

\paragraph{Real arithmetic for the matrix action.}

All matrix--vector multiplications associated with real nodes and conjugate
node pairs can be carried out in real arithmetic as in
section~\ref{sec:poly}. Real poles likewise require real shifted solves.
For a conjugate pole pair $\xi=\alpha+i\beta$, $\overline\xi$, $\beta\ne0$,
and a real right-hand side $u$, the reported implementation first solves
$(A-\xi I)g=u$ in complex arithmetic. The real-second-solve option obtains
the pair action by solving $(A-\alpha I)w=\operatorname{Re}g$ in real
arithmetic. When this option is not used, it solves $(A-\overline\xi I)w=g$ and takes the real
part. The real second solve is valid when $A-\alpha I$ is nonsingular: from $g=(A-\alpha I+i\beta I)w$ and real $w$, one obtains $\operatorname{Re}g=(A-\alpha I)w$. Proximity to singularity can make this option numerically unattractive, even though $A-\xi I$ is nonsingular.

\paragraph{An alternative pair identity.}
For completeness, the pair action also satisfies
\[
 w=\big((A-\xi I)(A-\overline\xi I)\big)^{-1}u
   =\frac{\operatorname{Im}((A-\xi I)^{-1}u)}{\beta}.
\]
Indeed, putting $B=A-\alpha I$ gives
$g=(B+i\beta I)(B^2+\beta^2I)^{-1}u$, and the vector
$(B^2+\beta^2I)^{-1}u$ is real. This identity does not require $B$ itself
to be invertible. It suggests using one complex solve per conjugate pair;
this optional simplification is not included in the reported timing
comparisons, and its rounding behavior, particularly for small $|\beta|$,
is separate from the exact-arithmetic identity.

\begin{theorem}\label{thm:rkspace}
The rational functions $q_0,\ldots,q_m$ defined in \eqref{defratbasis} form a
basis of $\{p/\ell_m:p\in\mathbb P_m\}$, where
$\ell_m(z)=\prod_{\ell=0}^{m-1}(z-\xi_\ell)$.
\end{theorem}
\begin{proof}
The nonsingular interpolation matrix in Theorem~\ref{thm:coefsratlinex}
implies that the $m+1$ functions are linearly independent. Since
$\ell_{d_j}$ divides $\ell_m$, each $q_j$ can be expressed with common
denominator $\ell_m$. Its resulting numerator has degree
$j+m-d_j\le m$ for $j\ge1$; for $q_0=1$ the numerator is $\ell_m$.
Thus all the functions belong to the stated $(m+1)$-dimensional space,
and linear independence proves the assertion.
\end{proof}

\begin{corollary}
Suppose that $\operatorname{span}\{b,Ab,\ldots,A^mb\}$ has dimension $m+1$.
Then $q_0(A)b,\ldots,q_m(A)b$ form a basis of
\[
 \mathcal Q_{m+1}(A,b)=\ell_m(A)^{-1}
       \operatorname{span}\{b,Ab,\ldots,A^mb\},
\]
and the final approximation $r_m(A)b=\sum_{j=0}^mc_jq_j(A)b$ lies in this
rational Krylov subspace.
\end{corollary}

For the same final $m+1$ nodes and $m$ prescribed poles, this interpolant is
identical in exact arithmetic to the interpolant evaluated by PAIN.
Indeed, multiplying by $\ell_m$ reduces both constructions to the unique
degree-at-most-$m$ polynomial interpolating $\ell_m(z_i)f(z_i)$ at the
$m+1$ nodes. Intermediate approximants can use different pole prefixes when conjugate pairs are kept together. The final interpolation function is unchanged, while its conjugation-preserving realization provides real coefficients and fixed large-vector storage.

If $D\supset W(A)$ and the poles are outside $D$, the scalar error transfers to the matrix action through \eqref{crouzeixbound}. The numerical study below separates this conditional error bound from the sampled scalar stopping tests used in practice.

\section{Numerical experiments}\label{sec:experiment}
The experiments compare paired and unpaired realizations, contrast their costs with existing Leja, Taylor, and Krylov methods, and examine the different inputs required by interpolation and two-pass methods. We report operation counts, action-stage timings where available, and differences from numerical reference vectors. The parameter choices and stopping tests are specified below; they are not identical forward-error targets.

\paragraph{Timing convention and interpretation.}\label{par:timing}
All reported runtimes exclude one-time preprocessing, including spectral estimation, construction of the scalar approximation region, and scalar AAA approximation and pole extraction. They measure the subsequent action computation, not the end-to-end cost of obtaining $f(A)b$ from $A$, $b$, and $f$ alone. Precomputed fixed-shift LU factors are reused outside their setup cost. Work carried out within a timed MATLAB backslash call for a changing shift remains part of that solve, and scalar updates interleaved with an action recurrence remain part of its runtime. This distinction is important when comparing fixed-pole and multi-pole methods. Preprocessing can be reused for several right-hand sides, but its cost must be added for a complete single-action comparison.

The regions used in the tests are selected from partial spectral information, sometimes with manual adjustment; they are not certified numerical-range enclosures. Equation~\eqref{crouzeixbound} gives the theoretical implication of an enclosing region, whereas the tables assess the observed performance of the chosen regions. Throughout, a reported relative reference error means $\|f_m-f_{\mathrm{ref}}\|_2/\|f_{\mathrm{ref}}\|_2$.

\subsection{Test problems}\label{subsec:tests}
We use 13 matrices from partial differential equations (PDEs), differential-algebraic equations (DAEs), and model order reduction (MOR), together with 13 network problems. Most are available from the SuiteSparse Matrix Collection \cite{kolodziej2019suitesparse}. The matrix \texttt{e40r0500} is from Matrix Market \cite{MatrixMarket}, and \texttt{aerofoil} comes from a finite-volume discretization in aerodynamics, supplied by Sebastian Timme \cite{Timme2012,RostamiXue18}. The matrix \texttt{Laplace2d} is based on the standard five-point finite-difference discretization on the unit square, with 192 uniform subintervals in each direction and without the inverse-square mesh-size factor. We abbreviate \texttt{windtunnel\_evap3d} as \texttt{windtun3d}. For the network problems, the tested operators are graph Laplacians: the first five graphs are undirected and the remaining eight are directed. Table~\ref{problemdscp} lists the matrix dimensions and recorded nonzero counts.

All experiments were performed in MATLAB R2023b Update 8 on a MacBook Pro with a 2.9 GHz dual-core Intel Core i5 processor, 16 GB of 1867 MHz DDR3 memory, and macOS Monterey 12.7.6.

\smallskip

\begin{table}
{\footnotesize
   \begin{center}
    \caption{Test-matrix dimensions and recorded nonzero counts.}\label{problemdscp}
  \begin{tabular}{|r|r|r|r|r|r|}
    \hline
    \multicolumn{3}{|c|}{PDE/DAE/MOR problems} & \multicolumn{3}{c|}{Network problems} \\ \hline
    name & size & nnz & name & size & nnz \\ \hline
    \texttt{tols4000} & 4000 & 9584 & \texttt{ca-HepTh} & 9877 & 61823 \\ 
    \texttt{GT01R} & 7980 & 430909 & \texttt{as-22july06} & 22963 & 119835 \\
    \texttt{graham1} & 9035 & 337542 & \texttt{wing} & 62032 & 305120 \\
    \texttt{inlet} & 11730 & 328323 & \texttt{delaunay\_n16} & 65536 & 458686 \\
    \texttt{epb1} & 14734 & 95053 & \texttt{usroads} & 129164 & 460034 \\ \cline{4-6}
    \texttt{aerofoil} & 16388 & 627983 & \texttt{Lederberg} & 8843 & 50438 \\
    \texttt{poli3} & 16955 & 37849 & \texttt{California} & 9664 & 25814 \\
    \texttt{e40r0500} & 17281 & 558361 & \texttt{p2pGnutella04} & 10879 & 50873 \\
    \texttt{ns3Da} & 20414 & 1680777 & \texttt{foldoc} & 13356 & 50873 \\
    \texttt{raefsky3} & 21200 &  1488768 & \texttt{p2pGnutella25} & 22687 & 77392 \\
    \texttt{af23560} & 23560 &  460598 & \texttt{cit-HepTh} & 27770 & 380538 \\
    \texttt{Laplace2d} & 36481 & 145160 & \texttt{enron} & 69244 & 343852 \\
    \texttt{windtun3d} & 40816 & 810218 & \texttt{NotreDame\_www} & 325729 & 1233616 \\ \hline
    \end{tabular}
    \end{center}
  }
\end{table}

\subsection{Polynomial approximations for entire functions}\label{subsec:polytest}
We compare Algorithm~\ref{alg:realpolylejamv} with real Leja interpolation for the exponential \cite{CaliariEtal04}, truncated Taylor methods for the exponential and trigonometric/hyperbolic functions \cite{AlMohy2011,AlMohy18}, restarted Arnoldi \cite{EiermannErnst06,AfanasjewEtal08}, and a two-pass sketched Krylov implementation. The comparisons separate scalar node selection, arithmetic realization, and large-vector storage.

For the polynomial tests we scale and shift the matrices as $A\leftarrow\alpha A+\beta I$ to obtain spectral regions of roughly unit size near $-1$. The manually chosen parameters are listed in Table~\ref{polyparameters}. A rectangle estimated from the spectrum is enclosed in the minimum-capacity ellipse used by the real Leja code. The same ellipse is the scalar approximation region for Algorithm~\ref{alg:realpolylejamv}. Its boundary is sampled at $2^{14}=16384$ equally spaced values of the angular parameter. The PDE/DAE/MOR examples include substantial imaginary spectral components and therefore exercise the conjugate-pair construction. The real Leja method uses nodes on the focal interval, while our method selects nodes on the sampled elliptic boundary.

\begin{table}[tp]
{\footnotesize
   \begin{center}
    \caption{Scaling and shift parameters for the polynomial tests, $A\leftarrow\alpha A+\beta I$.} \label{polyparameters}
  \begin{tabular}{|r|r|r|}
    \hline
    problem & $\alpha$ & $\beta$ \\ \hline
    \texttt{tols4000} & $2.1431\,10^{-4}$ & $-0.8$ \\
    \texttt{GT01R} & $5.0784\,10^{-5}$ & $-1$  \\
    \texttt{graham1} & $2.4096\,10^{-5}$ & $-0.02$ \\
    \texttt{inlet} & $1.0079$ & $-1$ \\
    \texttt{epb1} & $-10$ & $0$  \\
    \texttt{aerofoil} & $1.1939\,10^{-3}$ & $-0.5$ \\
    \texttt{poli3} & $-1.8349$ & $0.8$ \\
    \texttt{e40r0500} & $1.5267\,10^{-1}$ & $-2$ \\ 
    \texttt{ns3Da} & $-4.6512\,10^{-1}$ & $-1$ \\
    \texttt{raefsky3} & $2.5000\,10^{-6}$ & $0$ \\
    \texttt{af23560} & $3.8335\,10^{-3}$ & $-0.8$ \\
    \texttt{Laplace2d} & $-0.25$ & $0$ \\
    \texttt{windtun3d} & $-4.2194\,10^{-1}$ & $-0.002$ \\ \hline
    \end{tabular}
    \end{center}
  }
\end{table}

We first consider $f(z)=e^{hz}$, $\cos(hz)$, and $\sinh(hz)$ for $h=1,5,25$. The vector $b$ is generated by MATLAB's \texttt{randn} after \texttt{rng default}. Algorithm~\ref{alg:realpolylejamv} and real Leja interpolation stop when
\[
 \frac{\|f_{k+1}-f_k\|_2}{\|f_{k+1}\|_2}\leq2^{-30}\approx9.3\times10^{-10},
\]
where $f_k$ denotes the current action approximation at a completed update. The Taylor algorithms use their \texttt{single precision} tolerance setting, $2^{-24}\approx6\times10^{-8}$, with default parameter selection; this is a tolerance choice, not a change from double-precision arithmetic. These tests control different quantities, so we compare their observed reference errors as well as their operation counts.

\begin{table}[tp]
\centering
\setlength{\tabcolsep}{3pt}
{\footnotesize
\caption{Polynomial approximations to $\exp(hA)b$, $\cos(hA)b$, and $\sinh(hA)b$ for $h=1,5,25$. Each pair of rows gives matrix--vector product counts and relative reference errors. Real Leja interpolation is included for the exponential only.} \label{tab:polycomp}
\begin{tabular}{l|ccc|cc|cc}
\multirow{2}{*}{Problem} & \multicolumn{3}{c|}{$\exp(hA)b$} & \multicolumn{2}{c|}{$\cos(hA)b$} & \multicolumn{2}{c}{$\sinh(hA)b$} \\
& Alg.\,\ref{alg:realpolylejamv} & real Leja & Taylor & Alg.\,\ref{alg:realpolylejamv} & Taylor & Alg.\,\ref{alg:realpolylejamv} & Taylor \\ \hline
\multirow{4}{*}{\texttt{tols4000}} & $14/26/62$ & $7/13/32$ & $188/212/310$ & $13/22/50$ & $197/217/305$ & $14/24/62$ & $198/222/331$ \\
& $8\,10^{-14}$ & $1\,10^{-12}$ & $1\,10^{-10}$ & $2\,10^{-12}$ & $2\,10^{-14}$ & $8\,10^{-14}$ & $2\,10^{-14}$ \\ 
& $3\,10^{-13}$ & $2\,10^{-12}$ & $2\,10^{-10}$ & $7\,10^{-12}$ & $2\,10^{-11}$ & $9\,10^{-12}$ & $2\,10^{-11}$ \\ 
& $5\,10^{-11}$ & $1\,10^{-11}$ & $4\,10^{-10}$ & $2\,10^{-11}$ & $1\,10^{-09}$ & $1\,10^{-11}$ & $4\,10^{-11}$ \\ \hline
\multirow{4}{*}{\texttt{GT01R}} & $14/23/56$ & $7/12/27$ & $12/52/280$ & $14/22/44$ & $29/75/295$ & $14/24/58$ & $30/77/314$ \\
& $2\,10^{-14}$ & $2\,10^{-13}$ & $3\,10^{-11}$ & $2\,10^{-14}$ & $7\,10^{-15}$ & $2\,10^{-14}$ & $7\,10^{-15}$ \\ 
& $3\,10^{-12}$ & $2\,10^{-12}$ & $1\,10^{-10}$ & $2\,10^{-12}$ & $7\,10^{-13}$ & $6\,10^{-13}$ & $6\,10^{-13}$ \\ 
& $3\,10^{-11}$ & $2\,10^{-11}$ & $3\,10^{-10}$ & $1\,10^{-11}$ & $9\,10^{-11}$ & $9\,10^{-11}$ & $4\,10^{-10}$ \\ \hline
\multirow{4}{*}{\texttt{graham1}} & $13/21/45$ & $10/17/38$ & $14/37/132$ & $13/25/63$ & $29/75/225$ & $14/24/52$ & $30/70/210$ \\
& $8\,10^{-13}$ & $6\,10^{-12}$ & $5\,10^{-11}$ & $6\,10^{-13}$ & $6\,10^{-12}$ & $8\,10^{-13}$ & $1\,10^{-11}$ \\ 
& $2\,10^{-11}$ & $3\,10^{-10}$ & $2\,10^{-10}$ & $9\,10^{-13}$ & $2\,10^{-11}$ & $8\,10^{-11}$ & $5\,10^{-11}$ \\ 
& $5\,10^{-10}$ & $1\,10^{-09}$ & $2\,10^{-09}$ & $4\,10^{-12}$ & $2\,10^{-10}$ & $4\,10^{-10}$ & $7\,10^{-10}$ \\ \hline
\multirow{4}{*}{\texttt{inlet}} & $15/27/67$ & $12/25/66$ & $13/27/95$ & $15/27/73$ & $29/57/149$ & $14/28/62$ & $30/54/155$ \\
& $3\,10^{-13}$ & $1\,10^{-11}$ & $1\,10^{-11}$ & $4\,10^{-13}$ & $6\,10^{-14}$ & $5\,10^{-12}$ & $8\,10^{-14}$\\ 
& $1\,10^{-11}$ & $2\,10^{-11}$  & $6\,10^{-10}$ & $5\,10^{-12}$ & $2\,10^{-11}$ & $1\,10^{-12}$ & $6\,10^{-11}$ \\ 
& $5\,10^{-09}$ & $2\,10^{-10}$ & $1\,10^{-09}$ & $1\,10^{-09}$ & $2\,10^{-11}$ & $1\,10^{-10}$ & $1\,10^{-10}$ \\ \hline
\multirow{4}{*}{\texttt{epb1}} & $15/25/55$ & $11/20/45$ & $12/23/104$ & $15/25/65$ & $25/55/143$ & $14/24/54$ & $26/52/135$ \\
& $3\,10^{-14}$ & $8\,10^{-12}$ & $5\,10^{-11}$ & $1\,10^{-14}$ & $4\,10^{-12}$ & $3\,10^{-13}$ & $3\,10^{-12}$\\ 
& $1\,10^{-12}$ & $4\,10^{-11}$  & $2\,10^{-09}$ & $3\,10^{-11}$ & $2\,10^{-11}$ & $5\,10^{-12}$ & $1\,10^{-10}$ \\ 
& $1\,10^{-10}$ & $4\,10^{-10}$ & $1\,10^{-08}$ & $2\,10^{-12}$ & $4\,10^{-11}$ & $3\,10^{-11}$ & $3\,10^{-10}$ \\ \hline
\multirow{4}{*}{\texttt{aerofoil}} & $15/26/70$ & $8/13/31$ & $189/205/289$ & $14/28/63$ & $197/223/311$ & $16/28/70$ & $198/224/310$ \\
& $1\,10^{-13}$ & $3\,10^{-14}$ & $5\,10^{-11}$ & $3\,10^{-12}$ & $3\,10^{-13}$ & $3\,10^{-14}$ & $3\,10^{-13}$\\ 
& $3\,10^{-11}$ & $8\,10^{-12}$  & $2\,10^{-10}$ & $3\,10^{-12}$ & $4\,10^{-11}$ & $9\,10^{-12}$ & $3\,10^{-10}$ \\ 
& $5\,10^{-12}$ & $3\,10^{-11}$ & $5\,10^{-10}$ & $1\,10^{-10}$ & $4\,10^{-10}$ & $2\,10^{-10}$ & $2\,10^{-10}$ \\ \hline
\multirow{4}{*}{\texttt{poli3}} & $15/27/65$ & $13/23/61$ & $13/23/92$ & $15/27/69$ & $25/49/123$ & $16/28/68$ & $26/50/135$ \\
& $8\,10^{-13}$ & $7\,10^{-12}$ & $5\,10^{-12}$ & $1\,10^{-12}$ & $4\,10^{-12}$ & $7\,10^{-14}$ & $5\,10^{-12}$ \\ 
& $3\,10^{-12}$ & $8\,10^{-11}$ & $3\,10^{-10}$ & $2\,10^{-11}$ & $5\,10^{-12}$ & $4\,10^{-12}$ & $1\,10^{-11}$ \\ 
& $2\,10^{-10}$ & $2\,10^{-10}$ & $1\,10^{-07}$ & $3\,10^{-10}$ & $7\,10^{-11}$ & $4\,10^{-11}$ & $4\,10^{-10}$ \\ \hline
\multirow{4}{*}{\texttt{e40r0500}} & $15/29/75$ & $8/15/38$ & $13/38/145$ & $16/30/73$ & $29/65/187$ & $16/28/74$ & $30/61/155$ \\
& $1\,10^{-12}$ & $9\,10^{-13}$ & $7\,10^{-11}$ & $1\,10^{-13}$ & $4\,10^{-13}$ & $9\,10^{-14}$ & $6\,10^{-13}$ \\ 
& $7\,10^{-12}$ & $5\,10^{-12}$ & $6\,10^{-10}$ & $3\,10^{-12}$ & $3\,10^{-12}$ & $2\,10^{-11}$ & $9\,10^{-12}$ \\
& $4\,10^{-11}$ & $4\,10^{-11}$ & $5\,10^{-09}$ & $5\,10^{-11}$ & $3\,10^{-10}$ & $2\,10^{-11}$ & $2\,10^{-11}$ \\ \hline
\multirow{4}{*}{\texttt{ns3Da}} & $15/25/63$ & $7/13/31$ & $11/32/268$ & $14/26/58$ & $25/53/289$ & $14/28/68$ & $26/55/288$ \\
& $3\,10^{-14}$ & $1\,10^{-12}$ & $1\,10^{-10}$ & $4\,10^{-13}$ & $4\,10^{-13}$ & $4\,10^{-13}$ & $4\,10^{-13}$ \\
& $2\,10^{-11}$ & $1\,10^{-12}$ & $7\,10^{-10}$ & $1\,10^{-12}$ & $8\,10^{-12}$ & $5\,10^{-13}$ & $4\,10^{-11}$ \\
& $1\,10^{-10}$ & $2\,10^{-11}$ & $3\,10^{-09}$ & $3\,10^{-11}$ & $3\,10^{-10}$ & $2\,10^{-10}$ & $4\,10^{-10}$ \\ \hline
\multirow{4}{*}{\texttt{raefsky3}} & $13/20/37$ & $10/16/33$ & $15/39/139$ & $13/21/52$ & $29/77/221$ & $14/22/40$ & $30/70/195$ \\
& $1\,10^{-13}$ & $3\,10^{-12}$ & $3\,10^{-11}$ & $1\,10^{-13}$ & $3\,10^{-11}$ & $1\,10^{-14}$ & $3\,10^{-11}$ \\
& $2\,10^{-12}$ & $2\,10^{-10}$ & $6\,10^{-11}$ & $7\,10^{-12}$ & $7\,10^{-11}$ & $9\,10^{-13}$ & $5\,10^{-11}$ \\
& $4\,10^{-11}$ & $1\,10^{-10}$ & $1\,10^{-10}$ & $2\,10^{-12}$ & $4\,10^{-09}$ & $2\,10^{-11}$ & $1\,10^{-9}$ \\ \hline
\multirow{4}{*}{\texttt{af23560}} & $13/25/66$ & $7/12/29$ & $11/34/293$ & $15/26/58$ & $25/55/303$ & $14/26/70$ & $26/55/306$ \\
& $8\,10^{-12}$ & $5\,10^{-13}$ & $3\,10^{-11}$ & $3\,10^{-14}$ & $2\,10^{-13}$ & $3\,10^{-13}$ & $2\,10^{-13}$ \\
& $1\,10^{-11}$ & $2\,10^{-11}$ & $6\,10^{-11}$ & $2\,10^{-12}$ & $3\,10^{-12}$ & $1\,10^{-11}$ & $3\,10^{-11}$ \\
& $7\,10^{-11}$ & $6\,10^{-11}$ & $3\,10^{-10}$ & $2\,10^{-11}$ & $2\,10^{-10}$ & $3\,10^{-10}$ & $1\,10^{-09}$ \\ \hline
\multirow{4}{*}{\texttt{Laplace2d}} & $11/19/37$ & $10/18/35$ & $11/22/74$ & $11/21/50$ & $21/49/133$ & $11/19/37$ & $22/48/117$ \\
& $2\,10^{-11}$ & $3\,10^{-11}$ & $3\,10^{-10}$ & $3\,10^{-11}$ & $3\,10^{-10}$ & $2\,10^{-11}$ & $2\,10^{-10}$ \\
& $5\,10^{-11}$ & $3\,10^{-11}$ & $2\,10^{-09}$ & $1\,10^{-11}$ & $2\,10^{-10}$ & $5\,10^{-11}$ & $4\,10^{-10}$ \\
& $3\,10^{-10}$ & $3\,10^{-10}$ & $7\,10^{-09}$ & $8\,10^{-11}$ & $2\,10^{-08}$ & $8\,10^{-11}$ & $7\,10^{-10}$ \\ \hline
\multirow{4}{*}{\texttt{windtun3d}} & $13/23/47$ & $11/19/40$ & $315/344/430$ & $15/25/62$ & $293/331/503$ & $14/24/50$ & $294/333/503$ \\
& $2\,10^{-12}$ & $2\,10^{-11}$ & $2\,10^{-10}$ & $2\,10^{-14}$ & $3\,10^{-12}$ & $5\,10^{-14}$ & $2\,10^{-13}$ \\
& $6\,10^{-13}$ & $2\,10^{-11}$ & $9\,10^{-10}$ & $8\,10^{-13}$ & $5\,10^{-11}$ & $6\,10^{-12}$ & $1\,10^{-10}$ \\ 
& $8\,10^{-12}$ & $5\,10^{-10}$ & $9\,10^{-11}$ & $3\,10^{-12}$ & $9\,10^{-10}$ & $1\,10^{-10}$ & $9\,10^{-10}$ \\ \hline
\end{tabular}
}
\end{table}

\begin{table}[tp]
\centering
\setlength{\tabcolsep}{2.6pt}
{\scriptsize
\caption{Polynomial interpolation and two-pass sketched Krylov for $h=25$. The three rows per problem give matrix--vector product counts, action-stage times in seconds, and relative reference errors. The column $k_{\mathrm{orth}}$ records the empirically selected orthogonalization window, not a proven minimum. Preprocessing is excluded from the times.} \label{tab:polycompnew}
\begin{tabular}{l|ccc|ccc|ccc}
\multirow{1}{*}{Problem} & \multicolumn{3}{c|}{$\exp(hA)b$} & \multicolumn{3}{c|}{$\cos(hA)b$} & \multicolumn{3}{c}{$\sinh(hA)b$} \\ \hline
& Alg.\,\ref{alg:realpolylejamv} & two-pass & two-pass  & Alg.\,\ref{alg:realpolylejamv}  & two-pass & two-pass & Alg.\,\ref{alg:realpolylejamv} & two-pass &  two-pass \\
& (conj) & Krylov & window & (conj) & Krylov & window & (conj) & Krylov & window \\ 
& & & $k_{\mathrm{orth}}$ & & & $k_{\mathrm{orth}}$ & & & $k_{\mathrm{orth}}$ \\ \hline
\multirow{3}{*}{\texttt{raefsky3}} & $36$ & $62$ &  & $50$ & $86$ &  & $42$ & $64$ &  \\
			& $0.106$ & $0.335$ & $5$ &  $0.188$ & $0.606$ & $5$ & $0.148$ & $0.446$ & $5$ \\
			& $7\,10^{-11}$ & $1\,10^{-10}$ &  & $8\,10^{-11}$ & $2\,10^{-10}$ &  & $6\,10^{-12}$ & $6\,10^{-11}$ &  \\ \hline
\multirow{3}{*}{\texttt{ns3Da}} & $64$ & $100$ &  & $58$ & $88$ &  & $68$ & $104$ &  \\
			& $0.278$ & $0.710$ & $5$ &  $0.306$ & $0.777$ & $5$ & $0.342$ & $0.951$ & $5$ \\
			& $4\,10^{-11}$ & $9\,10^{-11}$ &  & $9\,10^{-11}$ & $7\,10^{-12}$ &  & $2\,10^{-11}$ & $5\,10^{-12}$ &  \\ \hline
\multirow{3}{*}{\texttt{af23560}} & $66$ & $92$ &  & $56$ & $80$ &  & $68$ & $90$ &  \\
			& $0.103$ & $0.297$ & $5$ &  $0.154$ & $0.445$ & $5$ & $0.094$ & $0.177$ & $10$ \\
			& $2\,10^{-11}$ & $8\,10^{-12}$ &  & $1\,10^{-11}$ & $4\,10^{-11}$ &  & $6\,10^{-11}$ & $4\,10^{-12}$ &  \\ \hline
\multirow{3}{*}{\texttt{graham1}} & $46$ & $70$ &  & $64$ & $78$ &  & $50$ & $52$ &  \\
			& $0.061$ & $0.133$ & $5$ &  $0.123$ & $0.342$ & $25$ & $0.094$ & $0.177$ & $10$ \\
			& $2\,10^{-10}$ & $6\,10^{-10}$ &  & $3\,10^{-11}$ & $7\,10^{-11}$ &  & $5\,10^{-11}$ & $2\,10^{-10}$ &  \\ \hline
\multirow{3}{*}{\texttt{GT01R}} & $56$ & $68$ &  & $44$ & $50$ &  & $58$ & $66$ &  \\
			& $0.079$ & $0.162$ & $20$ &  $0.089$ & $0.156$ & $20$ & $0.124$ & $0.244$ & $30$ \\
			& $3\,10^{-11}$ & $1\,10^{-12}$ &  & $4\,10^{-11}$ & $4\,10^{-12}$ &  & $3\,10^{-11}$ & $3\,10^{-11}$ &  \\ \hline
\multirow{3}{*}{\texttt{tols4000}} & $62$ & $102$ &  & $48$ & $80$ &  & $62$ & $96$ &  \\
			& $0.046$ & $0.131$ & $50$ & $0.087$ & $0.175$ & $35$ & $0.093$ & $0.264$ & $45$ \\
			& $1\,10^{-10}$ & $1\,10^{-09}$ &  & $2\,10^{-10}$ & $5\,10^{-09}$ &  & $4\,10^{-11}$ & $1\,10^{-08}$ &  \\ \hline
\end{tabular}
}
\end{table}

\begin{table}[tp]
\centering
\setlength{\tabcolsep}{2.6pt}
{\scriptsize
\caption{Polynomial interpolation with and without conjugate completion, real Leja interpolation for the exponential, and restarted Arnoldi (\texttt{funm\_kryl}, restart dimension 20), for $h=50$. Rows give matrix--vector product counts, action-stage times in seconds, and relative reference errors. Preprocessing is excluded.} \label{tab:polycomp2}
{
\begin{tabular}{l|cccc|ccc|ccc}
\multirow{1}{*}{Problem} & \multicolumn{4}{c|}{$\exp(hA)b$} & \multicolumn{3}{c|}{$\cos(hA)b$} & \multicolumn{3}{c}{$\sinh(hA)b$} \\ \hline
& Alg.\,\ref{alg:realpolylejamv} & Alg.\,\ref{alg:realpolylejamv} & real & \texttt{funm} & Alg.\,\ref{alg:realpolylejamv}  & Alg.\,\ref{alg:realpolylejamv} & \texttt{funm} & Alg.\,\ref{alg:realpolylejamv} & Alg.\,\ref{alg:realpolylejamv} & \texttt{funm} \\
& (conj) & (no-conj) & Leja & \texttt{kryl} & (conj) & (no-conj) & \texttt{kryl} & (conj) & (no-conj) & \texttt{kryl} \\ \hline
\multirow{3}{*}{\texttt{aerofoil}} & $112$ & $111$ & $49$ & $100$ & $104$ & $103$ & $80$ & $118$ & $116$ & $80$ \\
			& $0.246$ & $0.314$ & $0.222$ & $0.149$ & $0.237$ & $0.297$ & $0.127$ & $0.256$ & $0.341$ & $0.130$ \\
			& $2\,10^{-10}$ & $3\,10^{-10}$ & $1\,10^{-10}$ & $3\,10^{-14}$ & $7\,10^{-11}$ & $3\,10^{-10}$ & $2\,10^{-13}$ & $5\,10^{-11}$ & $7\,10^{-11}$ & $4\,10^{-10}$ \\ \hline
\multirow{3}{*}{\texttt{e40r0500}} & $120$ & $118$ & $60$ & $120$ & $112$ & $111$ & $100$ & $118$ & $115$ & $100$ \\
			& $0.250$ & $0.340$ & $0.238$ & $0.189$ & $0.252$ & $0.352$ & $0.156$ & $0.249$ & $0.322$ & $0.144$ \\
			& $3\,10^{-11}$ & $1\,10^{-10}$ & $6\,10^{-11}$ & $8\,10^{-15}$ & $5\,10^{-11}$ & $1\,10^{-10}$ & $4\,10^{-12}$ & $8\,10^{-12}$ & $4\,10^{-10}$ & $8\,10^{-12}$ \\ \hline
\multirow{3}{*}{\texttt{ns3Da}} & $104$ & $103$ & $50$ & $80$ & $90$ & $88$ & $80$ & $108$ & $107$ & $80$ \\
			& $0.487$ & $0.623$ & $0.373$ & $0.471$ & $0.416$ & $0.531$ & $0.330$ & $0.511$ & $0.658$ & $0.327$ \\
			& $9\,10^{-11}$ & $3\,10^{-10}$ & $4\,10^{-11}$ & $4\,10^{-10}$ & $1\,10^{-10}$ & $8\,10^{-10}$ & $5\,10^{-14}$ &  $1\,10^{-10}$ & $3\,10^{-10}$ & $3\,10^{-13}$ \\ \hline
\multirow{3}{*}{\texttt{raefsky3}} & $52$ & $50$ & $50$ & $60$ & $80$ & $79$ & $80$ & $54$ & $52$ & $60$ \\
			& $0.164$ & $0.198$ & $0.116$ & $0.158$ & $0.251$ & $0.339$ & $0.216$ & $0.189$ & $0.225$ & $0.154$ \\
			& $1\,10^{-11}$ & $2\,10^{-10}$ & $2\,10^{-10}$ & $7\,10^{-15}$ & $5\,10^{-10}$ & $6\,10^{-10}$ & $2\,10^{-12}$ & $1\,10^{-11}$ & $2\,10^{-10}$ & $7\,10^{-15}$ \\ \hline
\multirow{3}{*}{\texttt{windtun3d}} & $70$ & $68$ & $60$ & $60$ & $100$ & $94$ & $80$ & $74$ & $73$ & $60$ \\
			& $0.185$ & $0.259$ & ${0.101}$ & $0.158$ & $0.277$ & $0.357$ & $0.208$ & $0.202$ & $0.273$ & $0.143$ \\
			& $7\,10^{-11}$ & $2\,10^{-10}$ & $9\,10^{-10}$ & $6\,10^{-13}$ & $1\,10^{-11}$ & $2\,10^{-09}$ & $6\,10^{-12}$ & $9\,10^{-11}$ & $2\,10^{-10}$ &  $1\,10^{-13}$ \\ \hline
\multirow{3}{*}{\texttt{af23560}} & $110$ & $109$ & $46$ & $80$ & $88$ & $86$ & $60$ & $112$ & $111$ & $80$ \\
			& $0.232$ & $0.316$ & $0.211$ & $0.120$ & $0.189$ & $0.248$ & $0.096$ & $0.231$ & $0.308$ & $0.122$ \\
			& $5\,10^{-11}$ & $2\,10^{-10}$ & $2\,10^{-10}$ & $1\,10^{-15}$ & $1\,10^{-10}$ & $6\,10^{-10}$ & $1\,10^{-14}$ & $3\,10^{-10}$ & $3\,10^{-10}$ & $4\,10^{-15}$ \\ \hline
\end{tabular}
}
}
\end{table}

Table~\ref{tab:polycomp} reports matrix--vector product counts and relative reference errors. We do not report timings for this first comparison because the short runs are sensitive to scalar setup and background activity. Reference vectors were computed by real Leja interpolation with tolerance $10^{-14}$ for the exponential and by the Taylor method \cite{AlMohy18} with its \texttt{double precision} tolerance for the other functions. For example, Algorithm~\ref{alg:realpolylejamv} uses 14, 26, and 64 products for \texttt{tols4000} with $h=1,5,25$, respectively, and gives reference errors of approximately $8\times10^{-14}$, $3\times10^{-13}$, and $5\times10^{-11}$.

Under these settings, the Taylor methods use more matrix--vector products than Algorithm~\ref{alg:realpolylejamv}, particularly at $h=25$. The differences reach factors of 4--10 for \texttt{aerofoil} and \texttt{windtun3d}, and factors of 3--4 for several other examples. Norm-based parameter estimates can be conservative for nonnormal matrices. The ratios $\|A\|_1/\rho(A)$ are approximately $3.4\times10^3$, $1.6\times10^2$, and $2.3\times10^3$ for \texttt{tols4000}, \texttt{aerofoil}, and \texttt{windtun3d}, respectively, but less than 10 for the other matrices. Since the Taylor algorithms also use estimates involving powers of $A$, these ratios alone do not explain their parameter choices; the conclusion supported by the table is the lower product count of Leja interpolation for the stated settings.

For the exponential, real Leja interpolation often requires fewer products than Algorithm~\ref{alg:realpolylejamv}. On \texttt{tols4000}, \texttt{GT01R}, \texttt{aerofoil}, \texttt{e40r0500}, \texttt{ns3Da}, and \texttt{af23560}, the count is about half as large. This reflects the different node constructions rather than a general relation between a real node and a conjugate pair. The real Leja implementation also benefits from accurate exponential divided differences \cite{McCurdy1984}. Our block-interpolation construction is less specialized: it applies directly to other analytic conjugation-preserving functions, but its coefficient accuracy depends on the interpolation system.

To examine arithmetic and storage more directly, we compare the conjugation-preserving polynomial method, labeled ``conj'', with the ordinary complex Leja construction, labeled ``no-conj''. We also include real Leja interpolation for the exponential and the \texttt{funm\_kryl} implementation of restarted Arnoldi \cite{EiermannErnst06,AfanasjewEtal08}, with restart dimension 20 and no thick restart. The common update tolerance is $2^{-30}$. For this benchmark, \texttt{funm\_kryl} can additionally stop when its relative reference error is at most $2^{-30}$; this reference-based test is not available in an application with an unknown action. Timings are averages over 25 runs. Table~\ref{tab:polycompnew} uses $h=25$ for the two-pass comparison, and Table~\ref{tab:polycomp2} uses $h=50$ for the conjugation and restarted-Arnoldi comparison.

Table~\ref{tab:polycompnew} compares Algorithm~\ref{alg:realpolylejamv} with a two-pass Krylov implementation using randomized sketching and truncated Arnoldi. For each matrix and function, we manually tried orthogonalization-window sizes and retained the smallest satisfactory choice among those tested; the corresponding values are reported as $k_{\mathrm{orth}}$. This is the number of recent basis vectors used in local orthogonalization, as in the incomplete orthogonalization method \cite[section~6.4.2]{saad2003iterative}, not a proven minimum or the total storage of the algorithm.

The required window varies substantially. A window of 5 is sufficient in most of the \texttt{raefsky3}, \texttt{ns3Da}, and \texttt{af23560} tests; some \texttt{graham1} and \texttt{GT01R} tests use 20--30 vectors. For \texttt{tols4000}, the selected windows are 50, 35, and 45, close to full first-pass orthogonalization, and the resulting reference errors are still of order $10^{-9}$--$10^{-8}$. The interpolation method has a lower action-stage time in this table and uses no orthogonalization window. The two-pass implementation also performs basis regeneration and sketch-related updates, including repeated evaluation of the small projected function. Its important advantage is that it does not require our prescribed scalar approximation region. Thus the comparison identifies a practical exchange between region information and matrix-dependent orthogonalization, rather than an inherent failure of two-pass methods.

In Table~\ref{tab:polycomp2}, conjugate-completed and ordinary complex Leja interpolation have similar product counts for most matrices. The paired realization has lower action-stage times, demonstrating the benefit of real large-vector arithmetic in this comparison. The counts differ more for \texttt{af23560}: pairing requires more points for the exponential and hyperbolic sine, but fewer for the cosine. Thus the arithmetic advantage need not be accompanied by an identical interpolation sequence.

Restarted Arnoldi often gives smaller reference errors with comparable product counts and is faster in several entries of Table~\ref{tab:polycomp2}. Its implementation updates the action at restarts, whereas the interpolation code updates the action and sampled scalar data at every block. These online scalar operations are included in the interpolation timing and can be significant when a grid of 16384 points is comparable in size to the matrix dimension. The methods therefore differ in both update frequency and basis management. Algorithm~\ref{alg:realpolylejamv} retains fewer large vectors than even a restarted Arnoldi method with dimension 20, but this storage advantage does not imply a runtime advantage in every case.

\subsection{Rational approximations to \texorpdfstring{$\exp(-hA)b$}{exp(-hA)b}}
We next consider matrices with spectra in the right half-plane and a broad spectral spread, where polynomial approximation is less attractive. We scale and shift them as $A\leftarrow\alpha(A+\beta I)$ using Table~\ref{expparameters}. Some shifts move eigenvalues closer to the imaginary axis, widening the sector needed for a scalar approximation to the exponential \cite{Moret.Novati.2004}.

The scalar region in Figure~\ref{fig:gensector} is specified by $\gamma$, $\theta$, and $\rho$: the partial disk has center a distance $\gamma$ from the origin, angular parameter $\theta$, and radius $\rho$. When $\gamma=0$, the region is a sector with opening angle $2\theta$, symmetric about the real axis. We choose the parameters from 20 eigenvalues of largest modulus, 20 of smallest modulus, and 20 with largest imaginary-part magnitude, computed by MATLAB's \texttt{eigs} with tolerance $10^{-6}$. Some parameters were then adjusted manually to obtain reliable convergence, without an optimization of the final degree. Table~\ref{expparameters} records the choices. These are spectrally informed approximation regions, not enclosures certified by the partial eigenvalue calculations. Related approximation-based methods also use problem-dependent scalar parameters or regions \cite{Gutteletal20,AlMohyHigham2010A}.

\begin{figure}[ht!]
\centering
\caption{The boundary of the scalar approximation region used for approximating $\exp(-hA)b$ by rational interpolation. Its shape and size are chosen from partial spectral information.} \label{fig:gensector}
\begin{tikzpicture}[scale= 0.90, transform shape]
    
    \draw[->] (-1,0) -- (3,0) node[right] {Re}; 
    \draw[->] (0,-2) -- (0,2) node[above] {Im};

    \fill[gray!20] (0,0.75) -- (1.5,0.75) arc (0:45:1.5) -- cycle; 
    \draw[thick] (1.5,0.75) arc (0:45:1.5);

    \fill[gray!20] (0,-0.75) -- (1.5,-0.75) arc (0:-45:1.5) -- cycle; 
    \draw[thick] (1.5,-0.75) arc (0:-45:1.5);

    \draw[thick] (0,0.75) -- (0,-0.75);

    \draw[thick] (1.5,0.75) -- (1.5,-0.75);

    \draw[thick] (0,0.75) -- (1.06,1.81);

    \draw[thick] (0,-0.75) -- (1.06,-1.81);

    \draw[decorate, decoration={brace, amplitude=5pt}, thick] (0.1,0.75) -- (0.1,0)
        node[midway, right, xshift=5pt] {$\gamma$};

    \node at (0.5,1.0) {$\theta$};

    \draw[decorate, decoration={brace, amplitude=5pt}, thick] (0.0,0.75) -- (1.06,1.81) 
        node[midway, above, yshift=5pt, xshift=-5pt] {$\rho$};

    \node[below left] at (0,0) {$0$};
\end{tikzpicture}
\end{figure}

\begin{table}[tp]
{\footnotesize
   \begin{center}
   \caption{Matrix scaling $A\leftarrow\alpha(A+\beta I)$ and scalar-region parameters for the rational exponential tests ($h=1$); see Figure~\ref{fig:gensector}.} \label{expparameters}
  \begin{tabular}{|r|c|c|c|c|c|}
    \hline
    problem & $\alpha$ & $\beta$ & $\gamma$ & $\theta$ & $\rho$  \\ \hline
    \texttt{GT01R} & $-1$ & $0$ & $5\pi$ & $0.97\pi/2$ & $2.5\,10^{4}$ \\ 
    \texttt{graham1} & $-0.01$ & $- 4.85702\,10^2$ & $0$ & $0.8\pi/2$ & $10^3$ \\ 
    \texttt{epb1} & $2000$ & $- 4.85285\,10^{-5}$ & $0$ & $0.35\pi/2$ & $5\,10^2$ \\ 
    \texttt{aerofoil} & $-10$ & $1.03999\,10^{-2}$ & $3\pi$ & $0.81\pi/2$ & $1.3\,10^{4}$ \\ 
    \texttt{e40r0500} & $10$ & $2.7\,10^{-2}$ & $0$ & $0.65\pi/2$ & $1.4\,10^2$  \\ 
    \texttt{ns3Da} & $100$ & $9.91908\,10^{-1}$ & $2\pi$ & $0.9\pi/2$ & $3\,10^2$ \\ 
    \texttt{raefsky3} & $-1$ & $0$ & $1.5$ & $0$ & $10^6$ \\ 
    \texttt{af23560} & $-5$ & $2.7306\,10^{-1}$ & $0$ & $0.91\pi/2$ & $1.5\,10^{3}$ \\ 
    \texttt{Laplace2d} & $100$ & $-5\,10^{-4}$ & $0$ & $0$ & $8\,10^2$  \\ 
    \texttt{windtun3d} & $1000$ & $0$ & $5\pi$ & $0.5\pi/2$ & $6\,10^3$ \\ \hline 
    \end{tabular}
    \end{center}
  }
\end{table}

\begin{table}[tp]
\centering
\setlength{\tabcolsep}{3pt}
{\footnotesize
\caption{Rational approximations to $\exp(-A)b$. Rows give linear-solve counts, action-stage times in seconds, and relative reference errors. \texttt{conj} and \texttt{no-conj} denote paired and unpaired realizations. Preprocessing is excluded.} \label{tab:exp_exp}
\begin{tabular}{l|cccc|ccc}
&  & & & & AAA\,+ & AAA\,+ & AAA\,+ \\
Problem & shift (fixed  & shift-invert & Alg.\,\ref{alg:realpolylejasimv} & Alg.\,\ref{alg:realpolylejasimv} & Alg.\,\ref{alg:realratlejamv} & Alg.\,\ref{alg:realratlejamv} & PAIN \\
& real pole) & Arnoldi & (conj) & (no-conj)  & (conj) & (no-conj) & (no-conj)\\ \hline
\multirow{3}{*}{\texttt{GT01R}} & & $467$ & $451$ & $434$ & $156$ & $155$ & $155$  \\
			& $-250$ & $5.57$ & $4.01$ & $4.14$ & $30.44$ & $34.16$ & $38.99$ \\
			& & $1\,10^{-08}$ & $5\,10^{-08}$ & $5\,10^{-08}$ & $5\,10^{-08}$ & $5\,10^{-08}$ & $5\,10^{-08}$ \\ \hline
\multirow{3}{*}{\texttt{graham1}} & & $44$ & $80$ & $76$ & $28$ & $28$ & $28$  \\
			& $-80$ & $0.31$ & $0.49$ & $0.54$ & $3.72$ & $4.26$ & $4.44$ \\
			& & $4\,10^{-10}$ & $3\,10^{-11}$ & $6\,10^{-10}$ & $6\,10^{-11}$ & $1\,10^{-10}$ & $1\,10^{-10}$ \\ \hline
\multirow{3}{*}{\texttt{epb1}} & & $29$ & $34$ & $31$ & $14$ & $12$ & $12$  \\
			& $-35$ & $0.25$ & $0.26$ & $0.28$ & $1.07$ & $1.08$ & $1.16$ \\
			& & $8\,10^{-11}$ & $3\,10^{-11}$ & $2\,10^{-10}$ & $8\,10^{-10}$ & $6\,10^{-09}$ & $6\,10^{-09}$ \\ \hline
\multirow{3}{*}{\texttt{aerofoil}} & & $86$ & $103$ & $99$ & $34$ & $33$ & $33$  \\
			& $-110$ & $1.79$ & $1.83$ & $1.93$ & $18.52$ & $20.11$ & $20.74$ \\
			& & $8\,10^{-09}$ & $3\,10^{-09}$ & $2\,10^{-09}$ & $2\,10^{-09}$ & $6\,10^{-09}$ & $6\,10^{-09}$ \\ \hline
\multirow{3}{*}{\texttt{e40r0500}} & & $47$ & $52$ & $52$ & $20$ & $19$ & $19$  \\
			& $-100$ & $0.54$ & $0.54$ & $0.60$ & $4.25$ & $4.63$ & $4.71$ \\
			& & $9\,10^{-10}$ & $2\,10^{-09}$ & $4\,10^{-09}$ & $1\,10^{-09}$ & $2\,10^{-07}$ & $3\,10^{-08}$ \\ \hline
\multirow{3}{*}{\texttt{ns3Da}} & & $71$ & $129$ & $130$ & $53$ & $52$ & $52$  \\
			& $-135$ & $4.03$ & $5.41$ & $6.06$ & $134.18$ & $150.86$ & $153.09$ \\
			& & $4\,10^{-10}$ & $3\,10^{-10}$ & $5\,10^{-10}$ & $1\,10^{-10}$ & $7\,10^{-10}$ & $7\,10^{-10}$ \\ \hline
\multirow{3}{*}{\texttt{raefsky3}} & & $23$ & $26$ & $25$ & $14$ & $11$ & $11$  \\
			& $-30$ & $1.35$ & $1.31$ & $1.37$ & $11.49$ & $10.65$ & $10.54$ \\
			& & $4\,10^{-09}$ & $6\,10^{-09}$ & $2\,10^{-09}$ & $3\,10^{-11}$ & $1\,10^{-09}$ & $1\,10^{-09}$ \\ \hline
\multirow{3}{*}{\texttt{af23560}} & & $83$ & $142$ & $136$ & $56$ & $55$ & $55$  \\
			& $-150$ & $2.24$ & $2.88$ & $3.02$ & $45.16$ & $50.95$ & $52.65$ \\
			& & $9\,10^{-10}$ & $1\,10^{-09}$ & $4\,10^{-09}$ & $4\,10^{-11}$ & $3\,10^{-10}$ & $3\,10^{-10}$ \\ \hline
\multirow{3}{*}{\texttt{Laplace2d}} & & $27$ & $28$ & $28$ & $10$ & $10$ & $10$  \\
			& $-40$ & $0.46$ & $0.37$ & $0.37$ & $1.55$ & $1.83$ & $1.88$ \\
			& & $1\,10^{-11}$ & $5\,10^{-11}$ & $5\,10^{-11}$ & $6\,10^{-09}$ & $6\,10^{-09}$ & $6\,10^{-09}$ \\ \hline
\multirow{3}{*}{\texttt{windtun3d}} & & $49$ & $61$ & $57$ & $25$ & $25$ & $25$  \\
			& $-80$ & $9.83$ & $10.71$ & $11.06$ & $221.36$ & $242.09$ & $239.05$ \\
			& & $6\,10^{-06}$ & $9\,10^{-07}$ & $9\,10^{-06}$ & $6\,10^{-10}$ & $5\,10^{-10}$ & $5\,10^{-10}$ \\ \hline
\end{tabular}
}
\end{table}

For $h=1$, we compare shift-invert Arnoldi (the restricted-denominator method \cite{Moret.Novati.2004}), Algorithm~\ref{alg:realpolylejasimv} with and without conjugate completion, and multi-pole interpolation using AAA with Algorithm~\ref{alg:realratlejamv} or PAIN. Shift-invert Arnoldi and the single-pole interpolation stop when $\|f_{k+1}-f_k\|_2/\|f_{k+1}\|_2\leq\mathrm{tol}$. The AAA-based methods use the sampled scalar condition
\[
 \frac{\|f(Z)-r_k(Z)\|_\infty}{\|f(Z)\|_\infty}\leq\mathrm{tol},
\]
where $Z$ samples the prescribed region and $r_k$ is the AAA approximant. The two stopping tests are useful computational criteria, but neither equates their matrix-action errors.

We first run shift-invert Arnoldi with tolerance $10^{-13}$ to obtain a reference vector, then compare the methods at tolerance $2^{-26}\approx1.5\times10^{-8}$. The multi-pole methods use MATLAB's backslash for shifted solves. Nodes and poles are ordered by increasing modulus, with conjugate pairs kept together in the paired realization. The same ordering is used for $A^{-1/2}b$, while decreasing modulus is used for $-A\log(A)b$. Other tested orderings, including Leja orderings, often produced more ill-conditioned interpolation systems or less regular convergence. The selected ordering is an empirical choice for these functions and regions, not a universal prescription.

Table~\ref{tab:exp_exp} reports solve counts, action-stage times, and reference errors. The maximum step count $m$ determines the fixed pole $s=-m/h=-m$. This standard choice is motivated by sectorial shift-invert error bounds \cite[equations~(2.3) and (3.23)]{Moret.Novati.2004}; it is not asserted to minimize the error for each individual matrix. For \texttt{GT01R}, the fixed-pole methods use four substeps with at most 250 solves per substep; all other examples use one substep. The fixed-shift matrix $A-sI$ is factored once by \texttt{lu}, and its sparse factors are reused.

Shift-invert Arnoldi and Algorithm~\ref{alg:realpolylejasimv} are broadly comparable, with the interpolation method generally using more steps. Their interpolation points arise differently: Arnoldi adapts through the generated subspace, whereas the explicit interpolation uses the prescribed scalar region. Multi-pole AAA interpolation uses fewer solves in these tests but is often substantially slower. Changing shifts requires more expensive solve work than reusing one real sparse factorization. An iterative-solver implementation could change that balance, but was not tested here. The short-recurrence methods retain a fixed number of large vectors; shift-invert Arnoldi retains the orthonormal basis, in addition to the solver storage needed by all methods.

Conjugate completion reduces the reported action-stage times relative to the corresponding unpaired constructions. The reduction is modest for the fixed-pole method, where both versions reuse the same real sparse factors. For the multi-pole implementation, roughly one third of the conjugate pole pairs in these exponential tests use the real-second-solve option described in section~\ref{sec:lmrksm}. This can reduce complex solve work. Paired and unpaired AAA can also select different poles, so their timing difference measures the combined interpolation and arithmetic choices rather than a pure real-versus-complex arithmetic experiment. The comparable solve counts indicate similar convergence for these cases. Agreement between PAIN and the unpaired Newton realization is consistent with their exact-arithmetic equivalence for the same final data.

\subsection{Rational approximations to \texorpdfstring{$A^{-1/2}b$}{the inverse-square-root action}}\label{subsec:sqrtinvAb}
We now consider $f(z)=z^{-1/2}$ on the principal branch. The successful methods in these tests are AAA-based multi-pole interpolation and adaptive rational Arnoldi with real poles for Cauchy--Stieltjes functions \cite{Guttel.Knizhnerman.2013}. The single-pole interpolation failed to converge for most matrices with the repeated real poles tested.

The sample set must avoid the singularity at zero. Starting from the region in Figure~\ref{fig:gensector}, we either remove boundary samples within distance $\delta_1$ of zero or translate the samples to the right by $\delta_2$. Table~\ref{invsqrtparameters} gives these parameters and records whether the matrix scaling and shift from Table~\ref{expparameters} were applied. The matrix modifications move unsuitable eigenvalues away from the branch cut. The subsequent changes to the scalar sample set are approximation choices, not proofs of a numerical-range enclosure. Both adaptive rational Arnoldi and the multi-pole methods use MATLAB's backslash for the shifted solves.

\begin{table}[tp]
{\footnotesize
   \begin{center}
     \caption{Scalar-region modifications for $A^{-1/2}b$. The interval $[10^2,10^6]$ is also added to the original boundary for \texttt{raefsky3}.}\label{invsqrtparameters}
  \begin{tabular}{|r|c|c|c|c|c|c|}
    \hline
    problem & modified & $\gamma$ & $\theta$ & $\rho$ & $\delta_1$ & $\delta_2$ \\ \hline
    \texttt{GT01R} & No & $5\pi$ & $0.97\pi/2$ & $2.5\,10^4$ & $10^{-1}$ & \\
    \texttt{graham1} & Yes & $0$ & $0.8\pi/2$ & $10^3$ & & $10^{-2}$ \\
    \texttt{epb1} & No & $0$ & $0.35\pi/2$ & $10^1$ & $4\,10^{-5}$ & \\
    \texttt{aerofoil} & No & $3\pi$ & $0.81\pi/2$ & $1.2\,10^3$ & $10^{-2}$ & \\
    \texttt{e40r0500} & Yes & $0$ & $0.87\pi/2$ & $1.4\,10^2$ & $10^{-2}$ & \\ 
    \texttt{ns3Da} & Yes & $2\pi$ & $0.9\pi/2$ & $2.5\,10^2$ & & $10^{-2}$ \\
    \texttt{raefsky3} & No & 1.45 & 0 & $10^2$ &  & $6.5\,10^{-6}$ \\ 
    \texttt{af23560} & No & $0$ & $0.8\pi/2$ & $3\,10^2$ & $2\,10^{-1}$ & \\ \hline
    \end{tabular}
    \end{center}
  }
\end{table}

\begin{table}[tp]
\centering
\setlength{\tabcolsep}{3pt}
{\footnotesize
\caption{AAA-based interpolation and adaptive rational Arnoldi for $A^{-1/2}b$. Rows give linear-solve counts, action-stage times in seconds, and relative reference errors. Preprocessing is excluded; the reference uses the tighter-tolerance AAA construction described in the text.} \label{tab:exp_invsqrt}
\begin{tabular}{l|ccc|c}
 & AAA\,+ & AAA\,+ & AAA\,+ & adaptive \\
Problem & Alg.\,\ref{alg:realratlejamv} & Alg.\,\ref{alg:realratlejamv} & PAIN & RKSM  \\
 & (conj) & (no-conj)  & (no-conj) &  (real pole) \\ \hline
 \multirow{3}{*}{\texttt{GT01R}}  & $32$ & $30$ & $30$ & $75$ \\
						& $4.50$ & $6.87$ & $6.97$ & $13.90$ \\
						& $2\,10^{-8}$ & $2\,10^{-5}$ & $2\,10^{-5}$ & $1\,10^{-4}$ \\ \hline
 \multirow{3}{*}{\texttt{graham1}}  & $27$ & $27$ & $27$ & $23$ \\
						& $2.99$ & $3.88$ & $3.91$ & $2.54$ \\
						& $3\,10^{-10}$ & $3\,10^{-7}$ & $3\,10^{-7}$ & $4\,10^{-9}$ \\ \hline						
 \multirow{3}{*}{\texttt{epb1}}     & $23$ & $22$ & $22$ & $27$ \\
						& $1.58$ & $1.84$ & $1.82$ & $1.94$ \\
						& $3\,10^{-10}$ & $9\,10^{-8}$ & $9\,10^{-8}$ & $2\,10^{-8}$ \\ \hline						
 \multirow{3}{*}{\texttt{aerofoil}}  & $28$ & $27$ & $27$ & $44$ \\
						& $10.75$ & $15.94$ & $16.82$ & $18.36$ \\
						& $1\,10^{-8}$ & $3\,10^{-5}$ & $2\,10^{-5}$ & $2\,10^{-7}$ \\ \hline
 \multirow{3}{*}{\texttt{e40r0500}}  & $25$ & $25$ & $25$ & $75$ \\
						& $4.33$ & $6.14$ & $6.39$ & $15.86$ \\
						& $2\,10^{-9}$ & $2\,10^{-6}$ & $2\,10^{-6}$ & $3\,10^{-7}$ \\ \hline						
 \multirow{3}{*}{\texttt{ns3Da}}   & $28$ & $29$ & $29$ & $43$ \\
						& $43.92$ & $79.13$ & $80.79$ & $62.11$ \\
						& $2\,10^{-9}$ & $8\,10^{-6}$ & $8\,10^{-6}$ & $5\,10^{-8}$ \\ \hline						
 \multirow{3}{*}{\texttt{raefsky3}} & $31$ & $33$ & $33$ & $75$ \\
						& $17.22$ & $28.78$ & $29.24$ & $43.92$ \\
						& $1\,10^{-6}$ & $2\,10^{-5}$ & $2\,10^{-5}$ & $1\,10^{-7}$ \\ \hline
\multirow{3}{*}{\texttt{af23560}} & $21$ & $21$ & $21$ & $65$ \\
						& $12.38$ & $19.35$ & $19.67$ & $36.07$ \\
						& $1\,10^{-8}$ & $3\,10^{-7}$ & $3\,10^{-7}$ & $8\,10^{-8}$ \\ \hline					
\end{tabular}
}
\vspace{-0.16in}
\end{table}
A numerical reference vector is computed by conjugation-preserving AAA with tolerance $10^{-12}$ followed by Algorithm~\ref{alg:realratlejamv}. Because this uses the same construction at a tighter scalar tolerance, the resulting errors are consistency measurements and can miss errors common to both approximations. The comparison runs use AAA tolerance $10^{-1}2^{-26}\approx1.5\times10^{-9}$ and include adaptive rational Arnoldi with negative real poles. Table~\ref{tab:exp_invsqrt} shows that the AAA-based methods usually use fewer solves, while the paired and unpaired versions have similar solve counts.

The paired construction also gives smaller sampled residuals after conversion from barycentric to Newton form in these inverse-square-root tests, despite slightly larger interpolation-system condition numbers. A similar pattern occurs for $-z\log z$ below, but not for the exponential. This is an empirical observation about the selected nodes, computed poles, and coefficient solves; it is not a general accuracy consequence of conjugation symmetry.

A more direct source of the timing difference is the computed pole geometry. In the paired inverse-square-root tests, the poles lie on the negative real axis, so the matrix solves are real. Standard AAA instead produces many poles with small but non-negligible imaginary parts, requiring complex solves in the unpaired realization. Thus both pole selection and arithmetic contribute to the lower action-stage times. The same mechanism is observed in the directed-graph examples below.

\subsection{Rational approximations to \texorpdfstring{$-A\log(A)b$}{the negative A log(A) action}}\label{subsec:AlogAb}
Finally, we consider $g(A)b$ for $g(z)=-z\log z$ with $g(0)=0$. For a trace-normalized symmetric graph Laplacian, $\operatorname{trace}(g(A))$ is the von Neumann entropy, whose estimation involves matrix-function actions or quadratic forms \cite{Benzi.etal.23,Minello.etal.18,Frommer.etal.24}. Since a graph Laplacian is singular, $-A\log(A)$ is shorthand for $g(A)$: the value is zero on the semisimple zero eigenspace, and the principal branch is used on the nonzero spectral part. For a nonnegative-weight graph Laplacian, the semigroup $e^{-tA}$ is bounded and stochastic, excluding nontrivial Jordan blocks at zero. The directed graphs provide nonsymmetric action tests; the usual positive-semidefinite entropy interpretation applies to the undirected graphs. The holomorphic-neighborhood bound \eqref{crouzeixbound} is not invoked at the singular point zero.

We compare the multi-pole methods, adaptive rational Arnoldi, and the single-pole interpolation. The latter uses a problem-specific pole selected by trial and error among the tested values. Polynomial approximations were slow in this setting and are omitted. Multi-pole and adaptive rational Arnoldi solves use MATLAB's backslash; the single-pole method reuses the sparse triangular factors from one initial \texttt{lu} call, whose setup time is excluded.

For an undirected graph, the scalar grid consists of $2^{14}+1=16385$ logarithmically spaced points on $[\epsilon_{\mathrm{mach}}^2,\lambda_{\max}]$, where $\epsilon_{\mathrm{mach}}^2=2^{-104}\approx5\times10^{-32}$ and $\lambda_{\max}$ is the largest positive eigenvalue of the trace-normalized Laplacian. The lower endpoint samples close to, but not at, zero. For directed graphs, samples lie on the boundary of the region in Figure~\ref{fig:gensector}, with $\gamma=0$ and the parameters in Table~\ref{entropyparameters}.

\begin{table}[htp!]
{\footnotesize
   \begin{center}
  \caption{Scalar-region parameters for directed-graph Laplacians.} \label{entropyparameters}
  \begin{tabular}{|r|c|c|}
    \hline
    problem & $\theta$ & $\rho$  \\ \hline
    \texttt{Lederberg} & $0.1\pi/2$ & $2.8\,10^{-2}$ \\
    \texttt{California} & $0.25\pi/2$ & $1.2\,10^{-2}$ \\
    \texttt{p2pGnutella04} & $0.15\pi/2$ & $2.6\,10^{-3}$ \\
    \texttt{foldoc} & $0.1\pi/2$ & $2.4\,10^{-3}$ \\
    \texttt{p2pGnutella25} & $0.12\pi/2$ & $1.3\,10^{-3}$ \\
    \texttt{cit-HepTh} & $0.07\pi/2$ & $2.0\,10^{-3}$ \\
    \texttt{enron} & $0.1\pi/2$ & $5.2\,10^{-3}$ \\
    \texttt{NotreDame\_www} & $0.2\pi/2$ & $6.7\,10^{-5}$ \\ \hline
    \end{tabular}
    \end{center}
  }
\end{table}

\begin{table}[htp!]
\centering
\setlength{\tabcolsep}{3pt}
{\footnotesize
\caption{Rational approximations to $g(A)b$, where $g(z)=-z\log z$ and $g(0)=0$, for undirected-graph Laplacians. Rows give linear-solve counts, action-stage times in seconds, and relative reference errors. Preprocessing is excluded.} \label{tab:exp_entropy1}
\begin{tabular}{l|cc|cc|c}
& shift of & & AAA\,+ & AAA\,+ & adaptive \\
Problem & Alg.\,\ref{alg:realpolylejasimv} & Alg.\,\ref{alg:realpolylejasimv} & Alg.\,\ref{alg:realratlejamv} & PAIN & RKSM \\ \hline
\multirow{3}{*}{\texttt{ca-HepTh}} & & $247$ & $16$ & $16$ & $33$  \\
	& $-3\,10^{-4}$ & $0.89$ & $0.80$ & $0.87$ & $1.92$ \\
	& & $4\,10^{-7}$ & $1\,10^{-9}$ & $1\,10^{-9}$ & $3\,10^{-7}$ \\ \hline
\multirow{3}{*}{\texttt{as-22july06}} & & $93$ & $15$ & $15$ & $37$  \\
	& $-4\,10^{-4}$ & $0.56$ & $1.42$ & $1.43$ & $3.87$ \\
	& & $2\,10^{-8}$ & $1\,10^{-9}$ & $1\,10^{-9}$ & $1\,10^{-8}$ \\ \hline	
\multirow{3}{*}{\texttt{wing}} & & $159$ & $11$ & $11$ & $31$  \\
	& $-2\,10^{-4}$ & $9.90$ & $9.49$ & $8.81$ & $24.78$ \\
	& & $5\,10^{-8}$ & $6\,10^{-9}$ & $1\,10^{-9}$ & $1\,10^{-8}$ \\ \hline
\multirow{3}{*}{\texttt{delaunay\_n16}} & & $228$ & $12$ & $12$ & $34$  \\
	& $-2\,10^{-4}$ & $3.27$ & $2.54$ & $2.65$ & $7.64$ \\
	& & $3\,10^{-8}$ & $4\,10^{-9}$ & $4\,10^{-9}$ & $6\,10^{-8}$ \\ \hline
\multirow{3}{*}{\texttt{usroads}}   & & $200$ & $15$ & $15$ & $29$  \\
	& $-4\,10^{-4}$ & $3.08$ & $2.90$ & $2.81$ & $6.95$ \\
	& & $3\,10^{-7}$ & $2\,10^{-10}$ & $2\,10^{-10}$ & $5\,10^{-5}$ \\ \hline				  
\end{tabular}
}
\end{table}

\begin{table}[htp!]
\centering
\setlength{\tabcolsep}{3pt}
{\footnotesize
\caption{Rational approximations to $g(A)b$, where $g(z)=-z\log z$ and $g(0)=0$, for directed-graph Laplacians. Rows give linear-solve counts, action-stage times in seconds, and relative reference errors. Preprocessing is excluded.} \label{tab:exp_entropy2}
\begin{tabular}{l|ccc|ccc|c}
& & & & AAA\,+ & AAA\,+ & AAA\,+ & adaptive \\
& shift of & Alg.\,\ref{alg:realpolylejasimv} & Alg.\,\ref{alg:realpolylejasimv} & Alg.\,\ref{alg:realratlejamv} & Alg.\,\ref{alg:realratlejamv} & PAIN & RKSM \\
Problem & Alg.\,\ref{alg:realpolylejasimv} & (conj) & (no-conj) & (conj) & (no-conj) & (no-conj) & (real pole) \\ \hline
\multirow{3}{*}{\texttt{Lederberg}} & & $234$ & $267$ & $24$ & $30$ & $30$ & $33$ \\
 	& $-8\,10^{-5}$ & $0.87$ & $1.29$ & $1.19$ & $2.28$ & $2.34$ & $1.76$ \\
	& & $1\,10^{-8}$ & $1\,10^{-8}$ & $1\,10^{-8}$ & $3\,10^{-7}$ & $3\,10^{-7}$ & $3\,10^{-8}$\\ \hline			  
\multirow{3}{*}{\texttt{California}} & & $314$ & $312$ & $19$ & $17$ & $17$ & $35$ \\
	& $-4\,10^{-5}$ & $1.30$ & $1.77$ & $0.23$ & $0.28$ & $0.27$ & $0.84$ \\
	& & $9\,10^{-9}$ & $1\,10^{-8}$ & $3\,10^{-9}$ & $6\,10^{-9}$ & $6\,10^{-9}$ & $1\,10^{-5}$\\ \hline
\multirow{3}{*}{\texttt{p2p-Gnutella04}} & & $228$ & $209$ & $19$ & $18$ & $18$ & $25$ \\
	& $-1\,10^{-5}$ & $2.05$ & $1.97$ & $4.23$ & $7.09$ & $7.38$ & $5.47$ \\
	& & $3\,10^{-9}$ & $5\,10^{-8}$ & $3\,10^{-11}$ & $7\,10^{-8}$ & $7\,10^{-8}$ & $1\,10^{-8}$\\ \hline
\multirow{3}{*}{\texttt{foldoc}} & & $192$ & $216$ & $21$ & $30$ & $30$ & $26$ \\
	& $-8\,10^{-6}$ & $20.61$ & $23.77$ & $198.87$ & $501.37$ & $487.83$ & $244.51$ \\
	& & $4\,10^{-8}$ & $3\,10^{-8}$ & $4\,10^{-10}$ & $5\,10^{-9}$ & $5\,10^{-9}$ & $2\,10^{-6}$\\ \hline
\multirow{3}{*}{\texttt{p2p-Gnutella25}} & & $204$ & $200$ & $19$ & $17$ & $17$ & $24$ \\
	& $-4\,10^{-6}$ & $1.87$ & $2.10$ & $4.62$ & $7.40$ & $7.35$ & $6.13$ \\
	& & $5\,10^{-8}$ & $1\,10^{-7}$ & $3\,10^{-11}$ & $1\,10^{-9}$ & $1\,10^{-9}$ & $9\,10^{-10}$\\ \hline
\multirow{3}{*}{\texttt{cit-HepTh}} & & $256$ & $272$ & $17$ & $16$ & $16$ & $35$ \\
	& $-6\,10^{-6}$ & $12.11$ & $15.25$ & $106.61$ & $235.67$ & $231.45$ & $236.36$ \\
	& & $1\,10^{-8}$ & $2\,10^{-8}$ & $5\,10^{-9}$ & $1\,10^{-7}$ & $1\,10^{-7}$ & $4\,10^{-8}$\\ \hline							  	  \multirow{3}{*}{\texttt{enron}} & & $294$ & $295$ & $19$ & $17$ & $17$ & $39$ \\
	& $-1\,10^{-5}$ & $6.82$ & $7.59$ & $24.91$ & $41.19$ & $40.58$ & $46.25$ \\
	& & $3\,10^{-9}$ & $7\,10^{-9}$ & $7\,10^{-10}$ & $3\,10^{-9}$ & $3\,10^{-9}$ & $2\,10^{-8}$\\ \hline
\multirow{3}{*}{\texttt{NotreDame\_www}} & & $87^*$ & $64^*$ & $19$ & $18$ & $18$ & $43$ \\
	& $-3\,10^{-4}$ & $2.46^*$ & $2.43^*$ & $18.22$ & $22.01$ & $21.79$ & $58.16$ \\
	& & $5\,10^{-5}$ & $1\,10^{-4}$ & $9\,10^{-8}$ & $2\,10^{-8}$ & $2\,10^{-8}$ & $3\,10^{-6}$\\ \hline	
\end{tabular}
}
\end{table}

Tables~\ref{tab:exp_entropy1} and \ref{tab:exp_entropy2} give the undirected- and directed-graph results. In the undirected tests the sample set is real, both AAA variants generate negative real poles, and the paired and unpaired constructions coincide. AAA-based interpolation uses at most 16 solves for these graphs and at most 24 for the directed graphs, fewer than adaptive rational Arnoldi in the reported comparisons. Its action-stage times are also lower in these examples.

For the directed graphs, paired and unpaired constructions again have similar solve counts. Paired AAA produces real poles, whereas unpaired AAA produces poles with small nonzero imaginary parts, increasing the cost of the shifted solves. The observed accuracy advantage should be read together with this difference in computed interpolation data, rather than attributed to real arithmetic alone. The single-pole method fails on \texttt{NotreDame\_www}; its best intermediate reference error is of order $10^{-4}$--$10^{-5}$.

\subsection{Discussion}\label{subsec:discussion}
The clearest benefit of the proposed construction is the combination of fixed large-vector storage and conjugation-preserving arithmetic. It provides a common realization for polynomial, single-pole, and multi-pole interpolation without full-basis orthogonalization. The comparisons with unpaired realizations demonstrate useful reductions in action-stage time. Comparisons with Taylor, restarted Arnoldi, and adaptive rational Arnoldi also show that no one method dominates: function-specific approximation, the quality of the region, and shifted-solve costs can outweigh the recurrence savings.

The two-pass comparison identifies a different tradeoff. Its observed orthogonalization windows vary markedly across matrices and functions, whereas the interpolation recurrence has no such parameter. In exchange, interpolation needs useful scalar approximation information. This can be appropriate when problem-specific bounds are available or several right-hand sides share the preprocessing. All timing comparisons here exclude that preprocessing, and the laptop experiments do not measure parallel speedup.

A conservative numerical-range enclosure need not require an eigenvalue computation of $A$. For a real matrix, set
\[
 H=\tfrac12(A+A^T),\qquad K=\tfrac1{2i}(A-A^T).
\]
If $\ell\leq\lambda_{\min}(H)$, $\lambda_{\max}(H)\leq u$, and $\|K\|_2\leq v$, then $W(A)\subset[\ell,u]+i[-v,v]$. Entrywise Gershgorin bounds provide such scalars; tighter bounds require additional work. A conservative rectangle may nevertheless be too large, or intersect a singularity or branch cut of $f$. This elementary construction is an option for obtaining region information, not the procedure used in our spectrally informed tests.

The tests also reveal accuracy limitations in scalar interpolation and pole computation. The systems in \eqref{polylinex} and \eqref{ratlinex} can be ill-conditioned, and small barycentric residuals can deteriorate when converted to a Newton representation. In the original calculations, some nominally conjugate poles agreed to only three digits and were averaged; this changes the rational function and can affect its matrix action. Tighter AAA tolerances can also introduce Froissart doublets that are not fully removed by postprocessing. The exact-arithmetic results establish the interpolation space and real-coefficient structure, while the attainable floating-point accuracy depends on these scalar computations and on the shifted solves.

\section{Conclusion}\label{sec:conclusion}
We have developed conjugation-preserving polynomial and rational interpolation schemes for matrix-function actions with real matrices and vectors. Real nodes and nonreal conjugate pairs are handled in blocks, giving real interpolation coefficients and short recurrences with a fixed number of large vectors. The multi-pole construction uses AAA data and evaluates the same final interpolant as PAIN for identical admissible nodes and poles, while reorganizing the arithmetic to exploit conjugation.

The numerical comparisons demonstrate the practical benefit of this realization and clarify its scope. Conjugation-preserving realizations often reduce action-stage costs relative to the tested complex-arithmetic counterparts; factorization reuse can favor a single-pole method, while suitable multiple poles can substantially reduce solve counts. The two-pass comparison shows that orthogonalization-window choices can be matrix-dependent, whereas explicit interpolation instead requires a prescribed scalar region. These methods are therefore most attractive when large-vector storage matters and useful approximation data are already available or reusable. Their remaining practical limitations are the conditioning of scalar interpolation, the accuracy of computed poles, and the cost of shifted solves.

\section*{Conflict of interest}
This study does not have any conflicts to disclose.

\bibliography{ref_shortrec}

\end{document}